\newif\ifpreprint
\preprinttrue
\ifpreprint 
\documentclass[a4paper,natbib]{article}
\usepackage[affil-it]{authblk}
\usepackage[colorlinks,bookmarksopen,bookmarksnumbered,citecolor=red,urlcolor=red]{hyperref}
\usepackage[cm]{fullpage}
\usepackage{amsthm}
\newcommand{\tabsize}{}
\newcommand{\toprule}{\hline}
\newcommand{\midrule}{\hline}
\newcommand{\bottomrule}{\hline}
\else 
\documentclass[times]{fldauth}
\fi 

\newtheorem{theorem}{Theorem}

\usepackage{amssymb,amsmath,amsfonts}
\usepackage{graphicx}
\usepackage{bbm}
\usepackage{color}
\usepackage{algorithmic}
\usepackage{algorithm}
\usepackage{stmaryrd} 
\usepackage{cite}

\newcommand{\mat}[1]{\underline{#1}}
\newcommand{\gradient}{\operatorname{grad}}
\newcommand{\divergence}{\operatorname{div}}
\newcommand{\Id}{\operatorname{Id}}
\newcommand{\unittwo}{\ensuremath{\operatorname{Id}_{2\times 2}}}
\newcommand{\rhat}{\ensuremath{\hat{\vec{r}}}}
\newcommand{\bg}[1]{\ensuremath{#1_{\operatorname{ref}}}}
\newcommand{\dt}{\ensuremath{\mu \Delta t}}
\newcommand{\bnabla}{\boldsymbol{\nabla}}
\renewcommand{\vec}[1]{\ensuremath{\boldsymbol{#1}}}
\newcommand{\DDt}[1]{\frac{D#1}{Dt}}
\newcommand{\dr}{\partial_r}
\newcommand{\order}{\mathcal{O}}
\newcommand{\cellint}{\int_E dV\;}
\newcommand{\cellinthoriz}{\int_T d\rhat}
\newcommand{\cellintvert}{\int_{r_k}^{r_{k+1}} dr\;r^2\;}
\newcommand{\nbT}{\mathcal{N}(T)}
\newcommand{\Rearth}{R_{\operatorname{earth}}}
\newcommand{\Omegaearth}{\Omega_{\operatorname{earth}}}

\newcommand{\nhoriz}{n_{\hS}}
\newcommand{\TPMGfull}{\ensuremath{\operatorname{TPMG}^{\operatorname{(full)}}}}
\newcommand{\TPMGfac}{\ensuremath{\operatorname{TPMG}^{\operatorname{\otimes}}}}
\newcommand{\TPMGpartfac}{\ensuremath{\operatorname{TPMG}^{\operatorname{(partial)}}}}
\newcommand{\deltaA}{\delta A}
\newcommand{\Afac}{{A^{\otimes}}}
\newcommand{\Mfac}{{M^{\otimes}}}
\newcommand{\facNorm}[1]{||#1||_{\Afac}}
\newcommand{\Surface}{\mathcal{S}}
\newcommand{\hS}{{\Surface}}
\newcommand{\nablatwod}{\ensuremath{\bnabla_{\!\!\Surface}}}
\newcommand{\eHoriz}{{E^{\hS}(\phi)}}
\newcommand{\eVert}{{E^r(r)}}

\title{Efficient Multigrid Preconditioners for Atmospheric Flow Simulations at High Aspect Ratio}
\ifpreprint 
\author[1]{Andreas Dedner}
\author[*,2]{Eike M\"uller}
\author[2]{Robert Scheichl}
\affil[1]{Mathematics Institute, Zeeman Building, University of Warwick, Coventry CV4 7AL}
\affil[2]{Department of Mathematical Sciences, University of Bath, Bath BA2 7AY, United Kingdom}
\affil[*]{Email: \texttt{e.mueller@bath.ac.uk}}

\else 
\runningheads{A. Dedner et al.}{Multigrid Preconditioners for Atmospheric Flow at High Aspect Ratio}

\author{Andreas Dedner\affil{1}, Eike M\"{u}ller\affil{2}\corrauth\ and Robert Scheichl\affil{2}}

\address{\affilnum{1}Mathematics Institute, Zeeman Building, University of Warwick, Coventry CV4 7AL \break \affilnum{2}Department of Mathematical Sciences, University of Bath, Bath BA2 7AY, United Kingdom}

\corraddr{Department of Mathematical Sciences, University of Bath, Bath BA2 7AY, United Kingdom. E-mail: \texttt{e.mueller@bath.ac.uk}}

\cgsn{Natural Environment Research Council (NERC)}{NE/J005576/1, NE/K006754/1}

\fi 

\begin{document}
\ifpreprint 
\maketitle
\fi 
\begin{abstract}
Many problems in fluid modelling require the efficient solution of highly anisotropic elliptic partial differential equations (PDEs) in ``flat'' domains. For example, in numerical weather- and climate-prediction an elliptic PDE for the pressure correction has to be solved at every time step in a thin spherical shell representing the global atmosphere. This elliptic solve can be one of the computationally most demanding components in semi-implicit semi-Lagrangian time stepping methods which are very popular as they allow for larger model time steps and better overall performance. With increasing model resolution, algorithmically efficient and scalable algorithms are essential to run the code under tight operational time constraints. We discuss the theory and practical application of bespoke geometric multigrid preconditioners for equations of this type. The algorithms deal with the strong anisotropy in the vertical direction by using the tensor-product approach originally analysed by B\"{o}rm and Hiptmair [Numer. Algorithms, 26/3 (2001), pp. 219-234]. We extend the analysis to three dimensions under slightly weakened assumptions, and numerically demonstrate its efficiency for the solution of the elliptic PDE for the global pressure correction in atmospheric forecast models.
For this we compare the performance of different multigrid preconditioners on a tensor-product grid with a semi-structured and quasi-uniform horizontal mesh and a one dimensional vertical grid. The code is implemented in the Distributed and Unified Numerics Environment (DUNE), which provides an easy-to-use and scalable environment for algorithms operating on tensor-product grids. Parallel scalability of our solvers on up to 20,480 cores is demonstrated on the HECToR supercomputer.
\end{abstract}
\ifpreprint 
\textbf{keywords}: 
\newcommand{\keysep}{,}
\else 
\newcommand{\keysep}{;}
\keywords{
\fi
high aspect ratio flow\keysep\ multigrid\keysep\ elliptic PDEs\keysep\ atmospheric modelling\keysep\ parallelisation\keysep\ convergence analysis
\ifpreprint 
\\[1ex]
\else
}
\maketitle
\fi
\ifpreprint 
\textbf{AMS classifiers}:
65N55, 
65Y20, 
65F08, 
65Y05, 
35J57, 
86A10 
\fi 
\section{Introduction}
Highly efficient solvers for elliptic partial differential equations (PDEs) are required in many areas of fluid modelling, such as numerical weather- and climate- prediction (NWP), subsurface flow simulations \cite{Lacroix2003} and global ocean models \cite{Marshall1997,Fringer2006}. Often these equations need to be solved in ``flat'' domains with high aspect ratio, representing a subsurface aquifer or the Earth's atmosphere. In both cases the horizontal extent of the area of interest is much larger than the vertical size. For example, the Euler equations, which describe the large scale atmospheric flow, need to be integrated efficiently in the dynamical core of NWP codes like the Met Office Unified Model \cite{Davies05,Wood2013}. Many forecast centres such as the Met Office and European Centre for Medium-Range Weather Forecasts (ECMWF) use semi-implicit semi-Lagrangian (SISL) time stepping \cite{Kwizak1971,Robert1981} to advance the atmospheric fields forward in time because it allows for larger model time steps and thus better computational efficiency. However, this method requires the solution of a anisotropic elliptic PDE for the pressure correction in a thin spherical shell at every time step.
As the elliptic solve can account for a significant fraction of the total model runtime, it is important to use algorithmically efficient and parallel scalable algorithms.

Suitably preconditioned Krylov-subspace and multigrid methods (see e.g. \cite{Saad2003,Trottenberg2001}) have been shown to be highly efficient for the solution of elliptic PDEs encountered in numerical weather- and climate prediction (see \cite{Barros1989,Bates1990,Bowman1991,Adams1989,Adams1991,Adams1992,Skamarock1997,Thomas1997,Hess1997,Qaddouri2003,Davies05,Buckeridge2010,Buckeridge2011} and the comprehensive review in \cite{Mueller2013}). Multigrid methods are algorithmically optimal, i.e. the number of iterations required to solve a PDE to the accuracy of the discretisation error is independent of the grid resolution. However - as far as we are aware - multigrid algorithms are currently not widely implemented operationally in atmospheric models and one of the aims of this paper is to demonstrate that they can be used very successfully in fluid simulations at high aspect ratio. Whereas ``black-box'' algebraic multigrid (AMG) \cite{Brandt1984,Stueben1999} solvers such as the ones implemented in the DUNE-ISTL \cite{Blatt10} and Hypre \cite{FalgoutYang2002} libraries can be applied under very general circumstances on unstructured grids and automatically adapt to potential anisotropies, they suffer from additional setup costs and lead to larger matrix stencils on the coarse levels. On (semi-) structured grids which are typical in many atmospheric and oceanographic applications, geometric multigrid algorithms usually give much better performance, as they can be adapted to the structure of the problem by the developer. In contrast to AMG algorithms which explicitly store the matrix on each level, it is possible to use a matrix-free approach: instead of reading the matrix from memory, it is reconstructed on-the-fly from a small number of ``profiles''.  This leads to a more regular memory access pattern and significantly reduces the storage costs, in particular if these profiles can be factorised into a horizontal and vertical component. As the code is memory bandwidth limited this also has a direct impact on the performance of the solver.
Robust geometric multigrid methods adapt the smoother or coarse grid transfer operators to deal with very general anisotropies in the problem (see e.g. \cite{Bey1997,Hackbusch1997,Dendy1987,Oosterlee1997,Reisinger2004,
deZeeuw1990,Schaffer1998}).
However, this robustness comes at a price and these methods are often computationally expensive and difficult to parallelise.

In the problems we consider, the tensor-product structure of the underlying mesh and the grid-aligned anisotropy make it possible to use the much simpler but highly efficient tensor-product multigrid approach described for example in \cite{Dendy1992,Pflaum1997}: line-relaxation in the strongly coupled direction is combined with semi-coarsening in the other directions only. The implementation is straightforward: in addition to an obvious modification of the intergrid operators, every smoother application requires the solution of a tridiagonal system of size $n_r$ in each vertical column with $n_r$ grid cells. The tridiagonal solve requires $\order(n_r)$ operations and hence the total cost per iteration is still proportional to the total number of unknowns. The method is also inherently parallel as in atmospheric applications domain decomposition is typically in the horizontal direction only.

In \cite{BoermHiptmair1999} this method was analysed theoretically for equations with a strong vertical anisotropy on a two dimensional tensor-product grid. The authors show that optimal convergence of the tensor-product multigrid algorithm in two dimensions follows from the optimal convergence of the standard multigrid algorithm for a set of one-dimensional elliptic problems in the horizontal direction. While the original work in \cite{BoermHiptmair1999} applies in two dimensions, it has been extended to three dimensions in \cite{Buckeridge2010PhD} and the algorithm has recently been applied successfully to three dimensional problems in atmospheric modelling in \cite{Buckeridge2010,Buckeridge2011}.
Although the proof in \cite{BoermHiptmair1999} relies on the coefficients in the PDE to factorise exactly into horizontal-only and vertical-only contributions, we stress that this property is not required anywhere in the implementation of the multigrid algorithm. In practice we expect the algorithm to work well also for approximately factorising coefficients and under suitable assumptions we are able to also prove this rigorously. To demonstrate this numerically, we carry out experiments for the elliptic PDE arising from semi-implicit semi-Lagrangian time stepping in the dynamical core of atmospheric models such as the Met Office Unified Model \cite{Davies05,Wood2013}, where the coefficients only factorise approximately but the multigrid convergence is largely unaffected.
Alternatively, we also investigate approximate factorisations of the atmospheric profiles and then apply the tensor product multigrid algorithm to the resulting, perturbed pressure equation to precondition iterative methods for the original system, such as a simple Richardon iteration or BiCGStab \cite{vanderVorst1992}.
As the operator is usually ``well behaved'' in this direction (i.e. it is smooth and does not have large variations on small length scales), the multigrid algorithm will converge in a very small number of iterations.

An additional advantage of applying the multigrid method only to the perturbed problem with factorised profiles is the significant reduction in storage requirements for the matrix. As the algorithm is memory bound and the cost of a matrix application or a tridiagonal solve depends on the efficiency with which the matrix can be read from memory this leads to performance gains in the preconditioner: we find that the time per iteration can be reduced by around $20\%$, but this has to be balanced with a possibly worse convergence rate. Nevertheless, our numerical experiments show, that in some cases the factorised preconditioner can be faster overall. On novel manycore computer architectures, such as GPUs, where around 30-40 floating point operations can be carried out per global memory access, we expect the performance gains from this matrix-free tensor-product implementation to be even more dramatic. If the matrix is stored in tensor product format and the local stencil is calculated on-the-fly, the costs for the matrix construction can essentially be neglected compared to the cost of reading fields from memory. For example, carrying out a sparse matrix-vector product requires 1 global memory read and $n_s$ global writes per grid cell compared to $2n_s$ reads and 1 write if the $n_s$-point matrix stencil is stored explicitly - a speedup of almost a factor two. The benefits of this matrix-free implementation on GPUs has recently been shown in a similar context in \cite{Mueller2013a}.

In state-of-the-art global weather prediction models the horizontal resolution is of the order of tens of kilometres with the aim of reducing this to around one kilometre in the next decade (the number of vertical grid cells is typically around $100$). The resulting problems with $10^9-10^{11}$ degrees of freedom can only be solved on operational timescales if their scalability can be guaranteed on massively parallel computers. In addition to the sequential algorithmic performance we demonstrate the parallel scalability of our solvers on HECToR, the UK's national supercomputer which is hosted and managed by the Edinburgh Parallel Computing Centre (EPCC). We find that our solvers show very good weak scaling on up to 20,480 AMD Opteron cores and can solve a linear system with 11 billion unknowns in less than 5 seconds (reducing the residual by five orders of magnitude).

All our code is implemented in the Distributed and Unified Numerics Environment (DUNE) \cite{Bastian2008a,Bastian2008b}, which is an object oriented C++ library and provides easy to use interfaces to common parallel grid implementations such as ALUGrid \cite{Schupp1999,Dedner2004,Burri2005} and UG \cite{Bastian1997}. Due to the modular structure of the library and because we can rely on the underlying parallel grid implementations, the implementation of our solvers on tensor-product grids is straightforward. Throughout the code performance is guaranteed by using generic metaprogramming based on C++ templates.

\paragraph{Structure}
This paper is organised as follows. In Section \ref{sec:EllipticPDE} we describe the pressure correction equation arising in semi-implicit semi-Lagrangian time stepping in atmospheric models and discuss the discretisation of the resulting linear PDE with particular emphasis on the tensor-product structure of the grid. The theory of the tensor-product multigrid algorithm is reviewed in Section \ref{sec:TensorProductMultigrid} where we extend the analysis in \cite{BoermHiptmair1999} to three dimensions following \cite{Buckeridge2010PhD}. In this section we also prove the convergence of the preconditioned Richardson iteration for non-factorising profiles.
The grid structure and the discretisation of the equation as well as the implementation of our algorithms in the DUNE framework are described in Section \ref{sec:Implementation}. Numerical results for different test cases are presented together with parallel scaling tests in Section \ref{sec:NumericalResults}. We conclude and present ideas for future work in Section \ref{sec:Conclusions}.
\ifpreprint 
Some more technical aspects can be found in the appendices, in particular the finite-volume discretisation is described in detail in Appendix \ref{sec:Discretisation}.
\fi 
\section{Elliptic PDE for pressure correction in atmospheric models}
\label{sec:EllipticPDE}
The elliptic PDE which arises in semi-Lagrangian semi-implicit time stepping in atmospheric forecast models is derived for example in \cite{Wood2013} for the ENDGame dynamical core of the Unified Model. For simplicity (and in contrast to \cite{Wood2013}) the work in this paper is based on a finite volume discretisation of a continuous version of this PDE and in the following we outline the main steps in the construction of the corresponding linear algebraic problem.

The Euler equations describe large scale atmospheric flow as a set of coupled non-linear differential equations for the velocity $\vec{v}$, (Exner-) pressure $\pi$, potential temperature $\theta$ and density $\rho$.
\begin{equation}
  \begin{aligned}
 \DDt{\vec{v}} &= - c_p \theta \nabla \pi + \vec{R}_v \qquad \text{(Momentum equation)} \\[1ex]
  \DDt{\theta} &= R_\theta \qquad \text{(Thermodynamic equation)}\\[1ex]
  \DDt{\rho} &= - \rho \nabla\cdot v \qquad\text{(Mass conservation)}\\[1ex]
  \rho\theta &= \Gamma\pi^{\gamma} \qquad\text{(Equation of state)}
  \end{aligned}
  \label{eqn:EulerEquations}
\end{equation}
The $R$-terms describe external- and sub-gridscale- forcings such as gravity and unresolved convection. The constants $\Gamma$ and $\gamma$ are defined as
\begin{xalignat}{3}
  \Gamma &\equiv p_0/R_d, &
  \gamma &\equiv \frac{1-\kappa}{\kappa}, &
  \kappa &\equiv R_d/c_p,
\end{xalignat}
where $p_0$ is a reference pressure; $c_p$ and $R_d$ are the specific heat capacity and specific gas constant of dry air.
System \eqref{eqn:EulerEquations} can be written schematically for the state vector $\Phi = \left\{\vec{v},\pi,\theta,\rho\right\}$ as \begin{equation}
  \frac{D\Phi}{Dt}(\vec{x},t) = \mathcal{N}[\Phi(\vec{x},t)].
  \label{eqn:EulerSchematic}
\end{equation}
Advection is described in the semi-Lagrangian \cite{Robert1981} framework, i.e. material time derivatives $D\Phi/Dt$ are replaced by
\begin{equation}
  \frac{D\Phi}{Dt}(\vec{x},t) \mapsto \frac{\Phi^{(t+\Delta t)}(\vec{x})-\Phi^{(t)}(\vec{x}_D)}{\Delta t}
\label{eqn:semiLagrangian}
\end{equation}
where $\vec{x}_D$ is the departure point of a parcel of air at time $t$ which is advected to position $\vec{x}$ at time $t+\Delta t$. The right-hand-side of (\ref{eqn:EulerSchematic}) is treated semi-implicitly \cite{Kwizak1971}. Because of the small vertical grid spacing and the resulting large Courant number of vertical sound waves, vertical advection needs to be treated fully implicitly, but some of the other terms are evaluated at the previous time step and thus treated explicitly; we write $\mathcal{N}=\mathcal{N}^{(\text{impl.})}+\mathcal{N}^{(\text{expl.})}$. We use the $\theta$-method with off-centering parameter $\mu$ for implicit time stepping and replace
\begin{equation}
 \begin{aligned}
  \mathcal{N}[\Phi(\vec{x},t)] &=
  \mathcal{N}^{(\text{impl.})}[\Phi(\vec{x},t)] +
  \mathcal{N}^{(\text{expl.})}[\Phi(\vec{x},t)] \\
  &\mapsto
  \mu \mathcal{N}^{(\text{impl.})}[\Phi^{(t+\Delta t)}(\vec{x})] +
  (1-\mu) \mathcal{N}^{(\text{impl.})}[\Phi^{(t)}(\vec{x})]
+
  \mathcal{N}^{(\text{expl.})}[\Phi^{(t)}(\vec{x})]
  \end{aligned}
  \label{eqn:semiImplicit}
\end{equation}
and in the following we always assume that $\mu=\tfrac{1}{2}$ which corresponds to the scheme described in \cite{CrankNicolson1996}.
By eliminating the potential temperature, density and all velocities from the resulting equation, one (non-linear) equation for the pressure $\pi^{(t+\Delta t)}$ at the next time step can be obtained\footnote{Mathematically this is equivalent to forming the pressure Schur complement of the equation.}. To solve this equation via (inexact) Newton iteration, all fields are linearised around a suitable reference state (which can for example be the atmospheric fields at the previous time step) denoted by subscript ``ref''. To this end the pressure at the next time step is written as $\pi^{(t+\Delta t)}(\vec{x}) \equiv \pi(\vec{x}) = \bg{\pi}(\vec{x})+\pi'(\vec{x})$ with analogous expressions for $\theta^{(t+\Delta t)}$ and $\rho^{(t+\Delta t)}$; the reference velocities $\bg{\vec{v}}$ are assumed to be zero. It should, however, be kept in mind that the linearisation does not need to be ``exact'' as the non-linear equation can be solved with an inexact Newton iteration.
In particular, some terms can be moved to the right hand side which is equivalent to treating them explicitly or lagging them in the non-linear iteration. Naturally, there will be a tradeoff between faster convergence of the Newton iteration and the cost of the inversion of the linear operator; for example, in \cite{Wood2013} all couplings to non-direct neighbours, which can be large in the case of steep orography, are moved to the RHS to reduce the size of the stencil of the linear operator. While these considerations are relevant for the optimisation of the non-linear solve in a particular model, in this article we focus on the solution of the linear equation, which is the computationally most expensive component of the Newton iteration.

Once the Exner pressure $\pi^{(t+\Delta t)}$ has been calculated, the evaluation of the remaining atmospheric fields at the next time step is straightforward and does not require any additional (non-)linear solves.
In contrast to explicit time stepping methods the Courant number can be chosen significantly larger than 1, which makes semi-implicit semi-Lagrangian time stepping very popular in operational models. However, because of the short advective time scale and to ensure that large scale flow is described correctly, the Courant number is usually limited to around 10, i.e. the implicit time step size is no more than one order of magnitude larger than what would be allowed in an explicit method. To evaluate the overall performance of the method, the benefits of a larger time step would have to be balanced against the additional cost for the elliptic solve.
\subsection{Linear equation}
\label{sec:LinearEquation}
For ease of notation we simply write $\pi\equiv\pi^{(t+\Delta t)}$ in the
following and drop the time indices. Then the non linear equation for $\pi$
is of the form
\begin{equation}
  \mathcal{N}\left(\pi\right) = \mathcal{R}.
  \label{eqn:non_linear}
\end{equation}
To solve this equation iteratively we expand all fields around a reference state (which can, for example, be given by the fields at the previous time step) to obtain a linear operator $\mathcal{L}$. As discussed above, in practise some terms might be lagged in the non-linear iteration, i.e. moved to the right hand side of the linear equation. At each step $k$ of the nonlinear iteration we write $\pi_k=\bg{\pi}+\pi'_k$  for the approximate solution to (\ref{eqn:non_linear}) and update the pressure as follows:
\begin{equation*}
  \begin{aligned}
  \text{Solve}\quad\mathcal{L} \pi'_{k} &= \tilde{\mathcal{R}}_{k-1} := \left(\mathcal{R}-\mathcal{L}\bg{\pi}\right)-\delta \mathcal{N}\left(\pi_{k-1}\right)\quad\text{for}\;\pi'_k\quad
  \text{with}\quad \delta\mathcal{N} \equiv \mathcal{N}-\mathcal{L},\\
  \text{Update}\quad\pi_k &= \bg{\pi}+\pi'_k.
\end{aligned}
\end{equation*}
Every iteration requires the solution of a linear equation $\mathcal{L}\pi'_{k}=\tilde{\mathcal{R}}_{k-1}$ for the pressure correction $\pi'_k$, which we denote as $\pi'$ in the following. To construct the linear operator $\mathcal{L}$ we proceed as follows: starting from (\ref{eqn:EulerEquations}) the semi-Lagrangian framework in (\ref{eqn:semiLagrangian}) is used for horizontal advection and vertical advection is treated implicitly (to ensure that mass is exactly conserved, advection is treated implicitly in all three spatial dimensions in the mass equation). The right hand sides are expanded according to (\ref{eqn:semiImplicit}). We linearise around reference profiles $\bg{\theta}$, $\bg{\pi}$ and $\bg{\rho}$ which fulfil the equation of state $\bg{\rho}\bg{\theta}=\left(\bg{\pi}\right)^\gamma$, i.e. write $\theta=\bg{\theta}+\theta'$ etc. and assume that the velocity expansion is around zero $\bg{\vec{v}}=0$.
If we split up the velocity into a tangential- and vertical- component $\vec{v}=(\vec{v}_\hS,w)$ the time-discretised Euler equations in (\ref{eqn:EulerEquations}) finally become in a spherical geometry
\begin{eqnarray}
  \vec{v}_\hS &=& \vec{R}'_{u_\hS} - \dt c_p \frac{1}{r}\left(
    \bg{\theta}\nablatwod \pi' + (\nablatwod\bg{\pi})\theta'
  \right), \label{eqn:SISLmom_horiz}\\[1ex]
  w &=& R_w' - \dt c_p \left(\bg{\theta}\dr \pi' + (\dr \bg{\pi}) \theta'\right), \label{eqn:SISLmom_vert}\\[1ex]
  \theta' &=& R'_\theta - \dt (\dr\bg{\theta})w, \label{eqn:SISLthermodyn}\\[1ex]
  \rho' &=& R_\rho' - \dt\left(
    \frac{1}{r^2}\dr\left(r^2 \bg{\rho} w\right) + \frac{1}{r}\left(\nablatwod\cdot \left(\bg{\rho}\vec{v}_\hS\right)\right)
  \right),\label{eqn:SISLmass}\\[1ex]
  \pi' &=& \frac{\bg{\pi}}{\gamma}\left(\frac{\rho'}{\bg{\rho}} + \frac{\theta'}{\bg{\theta}}\right),
\label{eqn:SISLstate}
\end{eqnarray}
where $\partial_r=\langle\hat{\vec{n}},\bnabla\rangle$ is the normal component of the derivative and $\nablatwod=r\left(\bnabla - \langle\hat{\vec{n}},\bnabla\rangle\right)$ is the component tangential to a unit sphere $\hS$ with outer normal $\hat{\vec{n}}$.
Any terms that depend on the current time step are absorbed in the $R'$-terms. We then rewrite (\ref{eqn:SISLstate}) as a function of $\rho'$ and insert it together with (\ref{eqn:SISLmom_horiz}) into (\ref{eqn:SISLmass}) to obtain an equation with $w$, $\rho'$ and $\pi'$ only
\begin{equation}
 \begin{aligned}
  -\frac{\theta'}{\bg{\theta}} + \gamma \frac{\pi'}{\bg{\pi}} &=
  R''_\rho - \dt \frac{\dr(r^2\bg{\rho}w)}{r^2\bg{\rho}}\\&\qquad+\;\;(\dt)^2 c_p \frac{\nablatwod\cdot\left(\bg{\rho}\bg{\theta}(\nablatwod \pi')\right)+\nablatwod\cdot\left(\bg{\rho}(\nablatwod\bg{\pi})\theta'\right)}{r^2\bg{\rho}}.
  \end{aligned}
  \label{eqn:substitution1}
\end{equation}
By solving (\ref{eqn:SISLmom_vert}) and (\ref{eqn:SISLthermodyn}) for $w$ and $\theta'$ we obtain
\begin{xalignat}{2}
   w &= \bg{\Lambda} \left(f_2-\dt c_p \bg{\theta}(\dr \pi')\right), &
   \theta' &= \bg{\Lambda} \left(f_3+(\dt)^2 c_p \bg{\theta}(\dr \bg{\theta})(\dr \pi')\right)
\label{eqn:wtheta}
\end{xalignat}
where $\bg{\Lambda} = \left(1+(\dt)^2 \left(\bg{N}\right)^2\right)^{-1}$
arises from the implicit treatment of vertical advection and the (squared) vertical buoyancy (or Brunt-V\"{a}is\"{a}l\"{a}-) frequency is given by
\begin{equation}
  \left(\bg{N}\right)^2 = -c_p (\dr \bg{\pi}) (\dr \bg{\theta}) = g \frac{\dr \bg{\theta}}{\bg{\theta}}.
  \label{eqn:BuoyancyFrequency}
\end{equation}
The functions $f_2$ and $f_3$ only depend on the fields at the current time step.
We rescale the vertical coordinate $r$ by the radius of the earth $R_{\operatorname{earth}}$ and the potential temperature by a reference temperature $T_0$ at ground level to make it dimensionless. Finally, we multiply equation \eqref{eqn:substitution1} by $\bg{\rho}$ and denote the typical horizontal velocity by
\begin{equation*}
  c_h \equiv \sqrt{\gamma}c_s\qquad\text{where}\quad c_s=\sqrt{c_p T_0/\gamma}
\end{equation*} is the speed of sound in a parcel of air with temperature $T_0$. Furthermore we introduce the dimensionless quantity
\begin{equation}
\omega = \frac{c_h \dt}{R_{\operatorname{earth}}}.
\end{equation}
After eliminating $w$ and $\theta$ from (\ref{eqn:substitution1}) with the help of (\ref{eqn:wtheta}) we obtain a second order equation for the pressure correction $\pi'$:
\begin{equation}
  \begin{aligned}
  &-\omega^2 \left\{
    \bg{\Lambda}\bg{\rho}(\dr\bg{\theta})(\dr{\pi'}) + \frac{1}{r^2}\left[\dr\left(r^2\bg{\Lambda}\bg{\rho}\bg{\theta}(\dr\pi')\right)+ \nablatwod\cdot\left(\bg{\rho}\bg{\theta}(\nablatwod \pi')\right)\right]
  \right\} \\
  &-\omega^4 \frac{1}{r^2} \nablatwod\cdot\left(\bg{\Lambda}\bg{\rho}(\nablatwod \bg{\pi})(\dr \bg{\theta})(\dr \pi')\right)
  + \gamma \frac{\bg{\rho}}{\bg{\pi}}\pi' = RHS
  \end{aligned}\label{eqn:HelmholtzApp}
\end{equation}
The $\order(\omega^4)$ term arises due to the last term in (\ref{eqn:SISLmom_horiz}). In \cite{Wood2013} this term is not included in the linear operator since all terms which stem from reference profiles that do not depend exclusively on the vertical coordinate are neglected. To be consistent with this approach, the $\order(\omega^4)$ term is assumed to be moved to the right hand side of the linear equation in the following. The first two terms in the curly brackets are the sum of a vertical advection and a vertical diffusion term.

In contrast, in \cite{Wood2013}, the linear pressure correction equation is derived from the discretised Euler equations. However, it can be shown that (\ref{eqn:HelmholtzApp}) is identical to the continuum limit of equation (67) in \cite{Wood2013} if the latter is written down explicitly in spherical coordinates. Denoting the unknown pressure correction $\pi'$ by $u$, as is common in the mathematical literature, the elliptic operator can be written as
\begin{equation}
  \begin{aligned}
  \mathcal{L}u &= -\omega^2 \vec{\nabla} \cdot(\mat{\alpha}\vec{\nabla} u) - \omega^2 \vec{\xi}\cdot \vec{\nabla} u + \beta u\\[1ex]
  &=
  -\omega^2
  \begin{pmatrix}
    \dr
    ,&
    \frac{1}{r}\nablatwod
  \end{pmatrix}^T
  \begin{pmatrix}
    \alpha_r & 0 \\[1ex]
    0 & \alpha_\hS\unittwo
  \end{pmatrix}
  \begin{pmatrix}
    \dr\\[1ex]
    \frac{1}{r}\nablatwod
  \end{pmatrix}
  u
  - \omega^2
  \begin{pmatrix}
    \xi_r,& 0
  \end{pmatrix}^T
  \begin{pmatrix}
    \dr \\[1ex]
    \frac{1}{r}\nablatwod
  \end{pmatrix}
  u
  + \beta u
  \end{aligned}
\label{eqn:HelmholtzVector}
\end{equation}
where $\unittwo$ is the $2\times 2$ identity matrix. The equation is solved in a thin spherical shell, $\Omega=\Surface\times [1,1+H]$ and $H=D/\Rearth \ll 1$ is the ratio of the thickness of the atmosphere and the radius of the earth. The solution $u=u(r,\rhat)$ depends on the coordinates $r\in[1,1+H]$ and $\rhat\in\Surface$.
In contrast to global latitude-longitude grids, on quasi-uniform grids the ratio between the smallest and largest grid spacing is bounded. To ensure that the horizontal acoustic Courant number $\approx\omega/h$
 (where $h$ is the smallest grid spacing) remains unchanged as the horizontal resolution is increased, the time step size $\Delta t$ has to decrease linearly with $h$. A simple scaling argument shows that the vertical advection term is much smaller than the diffusion term at high resolution.

The functions $\alpha_r(r,\rhat)$, $\alpha_\hS(r,\rhat)$, $\xi_r(r,\rhat)$ and $\beta(r,\rhat)$ are referred to as ``profiles'' in the following and can be obtained from the background fields $\bg{\pi}$, $\bg{\theta}$ and $\bg{\rho}$ by comparing the elliptic operators in (\ref{eqn:HelmholtzApp}) and (\ref{eqn:HelmholtzVector}):
\begin{xalignat}{4}
    \alpha_r &= r^2 \bg{\Lambda} \bg{\rho}\bg{\theta}
                   \left( = r^2 \bg{\Lambda} \alpha_\hS\right),&
    \alpha_\hS &= \bg{\rho}\bg{\theta},&
    \xi_r &= \bg{\Lambda}\bg{\rho}(\dr\bg{\theta}),&
    \beta &= \gamma \frac{\bg{\rho}}{\bg{\pi}}.
  \label{eqn:profiles}
\end{xalignat}
\subsection{Iterative solvers}
After discretisation, the Helmholtz equation in (\ref{eqn:HelmholtzVector}) can be written as a large algebraic system of the form
\begin{equation}
  A\mathbf{u} = \mathbf{f}.
  \label{eqn:AlgebraicEquation}
\end{equation}
where the finite-dimensional field vector $\mathbf{u}$ represents the pressure correction in the entire atmosphere. If we assume that the horizontal resolution is around 1 kilometre and $\order(100)$ vertical grid levels are used, each atmospheric variable has $\order(10^{10})$ degrees of freedom.
Problems of this size can only be solved with highly efficient iterative solvers and on massively parallel computers. Current forecast models, such as the Met Office Unified Model, use suitable preconditioned Krylov subspace methods (see e.g. \cite{Saad2003} for an overview) such as BiCGStab  \cite{vanderVorst1992,Wesseling2001}.
Due to the flatness of the domain the equation is highly anisotropic: typical grid spacings in the horizontal direction are at the order of tens of kilometres, whereas the distance between vertical levels can be as small as a few metres close to the ground. While this anisotropy is partially compensated by the ratio  $\alpha_\hS/\alpha_r=r^{-2}\bg{\Lambda}^{-1}$ in (\ref{eqn:profiles}), it remains large in particular for small time steps $\Delta t$ for which $\bg{\Lambda}\rightarrow 1$ (recall that we chose units such that $r\approx 1$).

As discussed in the literature \cite{Skamarock1997,Thomas1997,Qaddouri2003,Davies05}, a highly efficient preconditioner for Krylov methods in this case is vertical line relaxation. This amounts to a block Jacobi or block SOR iteration where the degrees of freedom in one vertical column are relaxed simultaneously by solving a tridiagonal equation. However, (geometric) multigrid algorithms have also been considered by the atmospheric modelling community \cite{Barros1989,Bates1990,Bowman1991,Adams1989,Adams1991,Adams1992,Hess1997,Chen2001,Buckeridge2010,Buckeridge2011} and recently some of the authors have demonstrated their superior behaviour for a simplified model equation \cite{Mueller2013}.
\section{Tensor-product multigrid for anisotropic problems}\label{sec:TensorProductMultigrid}
Efficient al\-go\-rithms for the solution of anisotropic equations have been studied extensively in the multigrid literature. For general anisotropies in convection dominated problems, robust schemes have been designed by adapting the smoother (see e.g.
\cite{Bey1997,Hackbusch1997}) or the coarsening strategy and the restriction/prolongation operators (see e.g. \cite{deZeeuw1990,Schaffer1998}). For example, in \cite{Dendy1987,Oosterlee1997,Reisinger2004} alternating approximate plane- and line- smoothers are discussed. Alternatively, if algebraic multigrid (AMG) \cite{Brandt1984,Stueben1999} is used, the coarse grids and smoothers will automatically adapt to any anisotropies and the method can even be applied on unstructured grids. However, AMG has additional setup costs for the coarse grids and explicitly stores the coarse grid matrices. This has a significant impact on the performance in bandwidth-dominated applications. While these ``black-box'' approaches work well for very general problems and do not require anisotropies to be grid-aligned, they can be computationally expensive and difficult to parallelise.
The problem is simplified significantly in the case of grid-aligned anisotropies, which are typical in atmospheric- and ocean-modelling applications.
It has long been known that if the problem is anisotropic in one direction only, this can be dealt with effectively by either adapting the smoother or coarsening strategy (see e.g. \cite{Mulder1989,Mavriplis1997,Oosterlee1998} and also the discussion for simple anisotropic model problems in \cite{Trottenberg2001}).

Both methods can be combined as for example discussed in \cite{Dendy1992,Pflaum1997,BoermHiptmair1999} where the solution of two dimensional anisotropic problems with grid-aligned anisotropies is studied. By using line-relaxation in the $r$-direction together with semi-coarsening in the $x$-direction only, the multigrid solver is robust with respect to anisotropies in both the $x$- and $r$-direction as long as they are grid-aligned. In the following we will refer to multigrid algorithms which combine horizontal semi-coarsening with vertical line relaxation in the strongly coupled direction as \textit{tensor product multigrid} (TPMG) methods (both in 2D and in 3D).

In \cite{BoermHiptmair1999} the convergence of such a tensor-product multigrid solver for elliptic equations of the form
\begin{equation*}
  \mathcal{L}u = -\divergence \left(\mat{\alpha}\gradient u \right) =
  -\begin{pmatrix}
    \partial_r & \partial_x
  \end{pmatrix}
  \begin{pmatrix}
  \alpha_r(r,x) & 0 \\
  0 & \alpha_x(r,x)
  \end{pmatrix}
  \begin{pmatrix}
    \partial_r \\
    \partial_x
  \end{pmatrix}
  u(r,x) = f(r,x)
\end{equation*}
in a two dimensional domain $\Omega=[0,1]\times[0,1]$ is analysed under the assumption that the coefficients in the diagonal $2\times 2$ matrix can be factorised, i.e.
\begin{xalignat}{2}
  \alpha_r(r,x)&=\alpha_r^{r}(r)\alpha_r^{x}(x), &  \alpha_x(r,x)&=\alpha_x^{r}(r)\alpha_x^{x}(x).\label{eqn:alphaxrfac}
\end{xalignat}
The authors show that the tensor product multigrid algorithm applied to this problem converges uniformly provided the standard multigrid algorithm with point relaxation and uniform coarsening converges uniformly for one dimensional (horizontal) operators of the form
\begin{equation*}
  \mathcal{L}^x(\lambda_j) u^x(x) = -\partial_x \left(\alpha_r^x(x)\partial_x u^x(x)\right)+\lambda_j\;\alpha_x^x(x)u^x(x)
\end{equation*}
where the positive values $\lambda_j$ are the eigenvalues of the vertical Galerkin matrices. In particular, they analyse the strongly anisotropic case of $\alpha_r \gg \alpha_x$, which arises for example in the case of a polar grid on a disk with a (small) hole at the origin.

\subsection{Tensor-product multigrid preconditioners}
\label{sec:TPMGs}
Based on these observations, we propose two approaches for solving the pressure correction equation in (\ref{eqn:HelmholtzVector}). In both cases we use an iterative method such as a Richardson iteration or BiCGStab and precondition it with the tensor-product multigrid algorithm.
\paragraph{Tensor-product multigrid with full, non-factorising profiles (\TPMGfull)}
Often the profiles encountered in atmospheric flow simulations only factorise approximately. Although the theory in \cite{BoermHiptmair1999} applies only if the coefficient functions $\alpha_r(r,x)$ and $\alpha_x(r,x)$ can be written as the product of a vertical and a horizontal function as in (\ref{eqn:alphaxrfac}), this assumption is not used anywhere in the implementation.
Our numerical experiments demonstrate that good convergence can be achieved even in the non-factorising case where we use the full operator in the multigrid preconditioner.
\paragraph{Tensor-product multigrid with approximately factorised profiles (\TPMGfac)}
Given a set of profiles, we explicitly construct an approximate factorisation and use the resulting operator in the multigrid preconditioner; we apply a matrix-free approach, where the local stencil is reconstructed on-the-fly from the profiles. Depending on the quality of the factorisation, this may lead to a slight increase in the number of iterations of the underlying iterative solver. However, in terms of computational cost this increase is usually offset by a reduction in the amount of data that needs to be transferred from main memory. As the algorithm is memory bound, this will translate directly into performance gains. If the profiles factorise, for each horizontal cell and edge, $\order(1)$ entries which describe the horizontal coupling need to be stored. In the vertical direction four vectors of length $n_r$ are required \textit{for the entire grid}. Hence, in this case the matrix can be reconstructed from $\order(n_r)+\order(n_{\hS})$ values where $\nhoriz$ is the number of horizontal grid cells. This should be compared to $\order(n_r\times n_{\hS})$ data transfers for constructing the matrix in the \TPMGfull\ approach. We also prove formally in the following section that the  Richardson iteration with \TPMGfac\ preconditioner converges if the non-factorising part of the operator is small.
\subsection{Convergence of Tensor-product multigrid}
The convergence theory for factorising profiles in a spherical shell is a straightforward generalisation of the two dimensional case in \cite{BoermHiptmair1999} and is written down in detail in \cite{Buckeridge2010PhD} based on the multigrid convergence theory in \cite{Hackbusch1985}. In the following we outline the proof for 
the pressure correction in
(\ref{eqn:HelmholtzVector}). This is done in two steps: we first argue that if the profiles factorise and the advection term is dropped, the resulting symmetric positive definite equation can be solved efficiently with a tensor-product multigrid iteration.
We then show that if this factorising operator is used as a preconditioner for a Richardson iteration, the method converges also for the non-factorising equation provided the non-factorising contribution is sufficiently small.
\subsubsection{Factorising case}
\label{sec:ProofFactorising}
Consider the following PDE in the spherical shell $\Omega=\Surface\times [1,1+H]$:
\begin{equation}\label{eqn:pdetheory}
  \mathcal{L}^\otimes u=-\omega^2 \bnabla \cdot \left[\mat{\alpha}(r,\rhat) \bnabla u(r,\rhat)\right]+\beta(r,\rhat) u(r,\rhat) = f(r,\rhat)
\end{equation}
with $r\in [1,1+H]$, $\rhat \in \Surface$.
This should be compared to (\ref{eqn:HelmholtzVector}); for simplicity we do not consider the vertical advection term in this section, as it may in general destroy the positive definiteness of the problem. However, for high horizontal resolution this term is small and can be  treated as a perturbation.
We further assume that the $3\times 3$ matrix $\mat{\alpha}$ and the function $\beta$ have the following form that factorises into the product of a horizontal and a vertical function:
\begin{xalignat*}{2}
  \mat{\alpha} &=
  \begin{pmatrix}
    \alpha_r(r,\rhat) & 0 \\[1ex]
    0 & \mat{\alpha}_\Surface(r,\rhat)
  \end{pmatrix}
  =
  \begin{pmatrix}
    \alpha_r^r(r) \alpha_r^\Surface(\rhat) & 0 \\[1ex]
    0 & \alpha_\hS^r (r)\mat{\alpha}_\hS^\hS(\rhat)
  \end{pmatrix}, &
  \beta(r,\rhat) &= \beta^r(r)\beta^\hS(\rhat).
\end{xalignat*}
We also require that $\beta^\hS(\rhat)=\alpha_r^\hS(\rhat)$, which is satisfied for all factorisations that we use in our numerical experiments.
The $2\times 2$ matrix $\mat{\alpha}_\hS^\hS$ is required to be symmetric positive definite and we assume $\alpha_\hS^r,\alpha_r^r,\alpha_r^\hS,\beta^r>0$.

To discretise the problem, we choose finite element spaces $V^\hS$ over $\Surface$ and $V^r$ over $[1,1+H]$ and tensorise them to obtain the product space $V\equiv V^r\otimes V^\hS$ over $\Omega$. We write $n_r \equiv\dim V^r$ and $n_\hS \equiv \dim V^\hS$. For any two functions $u(r,\rhat)=u^r(r)u^\hS(\rhat)$, $v(r,\rhat)=v^r(r)v^\hS(\rhat)$ in $V$ the bilinear form $a:V \times V \to \mathbb{R}$ associated with the operator $\mathcal{L}^\otimes$ in \eqref{eqn:pdetheory} can be expressed in terms of the bilinear forms
\begin{equation}
  \begin{aligned}
  a^r(u^r,v^r) &\equiv \int_1^{1+H}
    \alpha_r^r(r)\partial_r u^r(r)\partial_r v^r(r)\;dr, \\
  b^r(u^r,v^r) &\equiv \int_1^{1+H} \beta^r(r)u^r(r)v^r(r)\;dr,\\
  m^r(u^r,v^r) &\equiv \int_1^{1+H} \alpha_\hS^r(r) u^r(r)v^r(r)\;dr,\\
  a^\hS(u^\hS,v^\hS) &\equiv \int_\Surface \left\langle
    \mat{\alpha}_\hS^\hS(\rhat)\nablatwod u^\hS(\rhat),\nablatwod v^\hS(\rhat)
  \right\rangle\;d\rhat \ \text{ and}\\
  m^\hS(u^\hS,v^\hS) &\equiv \int_\Surface \alpha_r^\hS(\rhat) u^\hS(\rhat) v^\hS(\rhat)\;d\rhat
  \end{aligned}
  \label{eqn:BilinearFormsTP}
\end{equation}
as
\begin{equation*}
a(u,v) = \left(\omega^2 a^r(u^r,v^r)+b^r(u^r,v^r)\right)m^\hS(u^\hS,v^\hS)+\omega^2 m^r(u^r,v^r)a^\hS(u^\hS,v^\hS).
\end{equation*}
Using the Kronecker product, the Galerkin-matrix representation $\Afac$ of the bilinear form $a(\cdot,\cdot)$ can then be expressed in terms of the Galerkin matrices of the bilinear forms in (\ref{eqn:BilinearFormsTP}), i.e.
\begin{equation*}
  \Afac = \left(\omega^2 A^r +B^r\right) \otimes M^\hS + \omega^2 M^r \otimes A^\hS.
\end{equation*}
Here $A^r,M^r,B^r\in \mathbb{R}^{n_r \times n_r}$ correspond to $a^r(\cdot,\cdot)$, $m^r(\cdot,\cdot)$ and $b^r(\cdot,\cdot)$ respectively and describe the vertical derivative- and mass- matrices. Analogously the derivative and mass matrix in the horizontal direction are described by $A^\hS,M^\hS\in \mathbb{R}^{n_\hS \times n_\hS}$, which correspond to $a^\hS(\cdot,\cdot)$ and $m^\hS(\cdot,\cdot)$.

To use the tensor-product multigrid approach, we further assume that there is a nested sequence
\begin{equation}
  V_0^\hS\subset V_1^\hS \subset \dots \subset V_L^\hS \equiv V^\hS
\end{equation}
of finite element spaces over $\Surface$,
where the subscript $\ell$ denotes the multigrid level; for the icosahedral and cubed sphere grid this hierarchy naturally exists.
We then use $V_\ell \equiv V^r \otimes V_\ell^\hS$ to discretise the full three dimensional problem on the multigrid level $\ell$, i.e. we do not coarsen in the vertical direction. The line smoother then corresponds to collectively relaxing all degrees of freedom in each of the $n_r$-dimensional subspaces $\operatorname{span}\left\{\left\{\psi_{\ell,k} \right\}\otimes V^r\right\}$ where $\psi_{\ell,k}$ are the nodal basis functions on level $\ell$.

The two-dimensional prolongation $P^\hS_\ell: V_\ell^\hS\rightarrow V_{\ell+1}^\hS$ and restriction $R^\hS_\ell\equiv \left(P^\hS_\ell\right)^T: V_{\ell+1}^\hS\rightarrow V_\ell^\hS$ naturally induce intergrid transfer operators between the three dimensional spaces $V_\ell$ and $V_{\ell+1}$ by $P_\ell \equiv \Id \otimes P_\ell^\hS$, $R_\ell \equiv \Id \otimes R_\ell^\hS$. On each multigrid level the matrix $A^\otimes_{\ell}$ can be constructed recursively using the Galerkin product $A^\otimes_{\ell}\equiv R_\ell A^\otimes_{\ell+1}P_\ell$ and it is easy to see that $A^\otimes_{\ell}$ and the (block-)smoother $W^\otimes_{\ell}$ can be written as
\begin{equation*}
\begin{aligned}
    A^\otimes_{\ell} &= \left(\omega^2 A^r +B^r \right)\otimes M^\hS_\ell + \omega^2 M^r \otimes A^\hS_\ell,\\
    W^\otimes_{\ell} &= \left(\omega^2 A^r +B^r \right)\otimes W^{M,\hS}_\ell + \omega^2 M^r \otimes W^{A,\hS}_\ell.
\end{aligned}
\end{equation*}
In the case of weighted block-Jacobi relaxation, for example, the matrices $W^{M,\hS}_\ell$ and $W^{A,\hS}_\ell$ are the weighted diagonals of $M^\hS_\ell$ and $A^\hS_\ell$. One V-cycle of the tensor product multigrid algorithm can now be written down compactly as follows.
\begin{algorithm}[H]
 \caption{\textbf{Tensor Product Multigrid V-cycle.} Input: system matrix $A^\otimes_{\ell}$, RHS $\mathbf{f}_\ell$, initial guess $\mathbf{u}_\ell$}
 \begin{algorithmic}[1]
\STATE{\textit{Pre-Smoothing:} $\nu^\text{pre}_\ell$ steps of $\mathbf{u}_\ell \ \to \ (W^\otimes_\ell)^{-1} \mathbf{f}_\ell + \big(I - (W^\otimes_\ell)^{-1} A^\otimes_{\ell}\big)\mathbf{u}_\ell$}
\IF{$\ell > 0$}
\STATE{\textit{Residual Calculation:} $\mathbf{r}_\ell = \mathbf{f}_\ell - A^\otimes_{\ell}\mathbf{u}_\ell$}
\STATE{\textit{Recursion:} Apply algorithm with $A^\otimes_{\ell-1}$, $\mathbf{f}_{\ell-1} = R_{\ell-1} \mathbf{r}_\ell$, $\mathbf{u}_{\ell-1} = 0$}
\STATE{\textit{Coarse Grid Correction:} $\mathbf{u}_\ell \ \to \ \mathbf{u}_\ell + P_{\ell-1} \mathbf{u}_{\ell-1}$}
\ENDIF
\STATE{\textit{Post-Smoothing:} $\nu^\text{post}_\ell$ steps of $\mathbf{u}_\ell \ \to \ (W^\otimes_\ell)^{-1} \mathbf{f}_\ell + \big(I - (W^\otimes_\ell)^{-1} A^\otimes_{\ell}\big)\mathbf{u}_\ell$}
 \end{algorithmic}
\end{algorithm}
On the finest level $L$, this V-cycle is applied to the right hand side $\mathbf{f}_L = \mathbf{f}$ of the original problem until the residual error is reduced below a certain tolerance. We typically choose the numbers of smoothing steps to be $\nu_\ell^\text{pre} = 2 $ and $\nu_\ell^\text{post} = 2$, for $\ell >0$, and $\nu_0^\text{pre} + \nu_0^\text{post} = 1$ on the coarsest grid. To simply apply a few steps of the smoother on the coarsest grid is sufficient because the CFL condition ensures that the system matrix $A^\otimes_{0}$ on the coarsest grid is dominated by the mass matrix term  $B^r_\ell \otimes M^\hS_\ell$ and thus well-conditioned.
\paragraph{Reduction of the theory to two dimensions}
The crucial idea in \cite{BoermHiptmair1999} is now that it is possible to construct a set of $n_r$ invariant $n_\hS$-dimensional subspaces such that the convergence of the tensor product multigrid method for the problem in $\Omega \subset \mathbb{R}^3$ can be analysed by independently studying the convergence of a standard multigrid algorithm in each of these subspaces over $\hS \subset \mathbb{R}^2$. This can be seen as follows: because both $A^r$ and $M^r$ are positive definite, there exists an eigenbasis $\mathbf{e}_j^r$, $j=1,\dots,n_r$, of $V^r$ and a corresponding set of strictly positive eigenvalues $\lambda_j$ such that
\begin{xalignat}{2}
  \left(\omega^2 A^r+B^r\right)\mathbf{e}_j^r &= \lambda_j M^r\mathbf{e}_j^r, &
  \left\langle M^r\mathbf{e}_j^r,\mathbf{e}_k^r\right\rangle = \delta_{j,k} \qquad\text{for all}\quad j,k\in\{1,\dots,N\}.
  \label{eqn:EigenvalueEquation}
\end{xalignat}
It follows from simple identities for the inner product on tensor product spaces that
\begin{equation*}
  \begin{aligned}
  \left\langle
    \mathbf{e}_k^r \otimes \mathbf{u}^\hS, A^\otimes_{\ell}(\mathbf{e}^r_j\otimes \mathbf{v}^\hS)
  \right\rangle
  &=
  \left\langle
    \mathbf{e}_k^r \otimes \mathbf{u}^\hS, \left(\omega^2 A^r+B^r\right) \mathbf{e}^r_j\otimes M_\ell^\hS \mathbf{v}^\hS + M^r \mathbf{e}_j^r \otimes A_\ell^\hS \mathbf{v}^\hS
  \right\rangle\\
  &= \delta_{jk} \left\langle
    \mathbf{u}^\hS,\left(\lambda_j M_\ell^\hS+A_\ell^\hS\right) \mathbf{v}^\hS
  \right\rangle
  \end{aligned}
\end{equation*}
and so the subspaces spanned by the different $\mathbf{e}^r_j$ are $A^\otimes_{\ell}$-orthogonal, with a similar property for the smoother matrix $W^\otimes_\ell$. As we do not coarsen in the vertical direction, the intergrid operators $P_\ell$ and $R_\ell$ do not mix different subspaces. For each $j$ the space $\operatorname{span}\{\mathbf{e}_j^r\}\otimes V_\ell^\hS$ is trivially isomorphic to $V_\ell^\hS$ and each of the $n_r$ independent subspaces corresponds to a two dimensional problem on $\Surface$ with the following matrix representation of the linear operator and smoother:
\begin{xalignat*}{2}
  A^\hS_{\ell,j} &\equiv \omega^2 A_\ell^\hS + \lambda_j M_\ell^\hS, &
  W^\hS_{\ell,j} &\equiv \omega^2 W^{A,\hS}_\ell + \lambda_j W^{M,\hS}_\ell.
\end{xalignat*}
In particular, $A^\hS_{\ell,j}$ is the Galerkin matrix which is obtained from discretising the bilinear form
$\omega^2 a^\hS({u}^\hS,{v}^\hS) + \lambda_j m^\hS({u}^\hS,{v}^\hS)$ on $V_\ell^\hS$. This bilinear form is the weak formulation of the following two dimensional operator:
\begin{equation}
  \mathcal{L}^\hS_j {u}^\hS(\rhat) = -\omega^2 \nablatwod \cdot\left(\mat{\alpha}_\hS^\hS(\rhat)\nablatwod {u}^\hS(\rhat)\right)+\lambda_j\alpha_r^\hS(\rhat){u}^\hS(\rhat)
\label{eqn:Cont2dOperator}
\end{equation}

\paragraph{Convergence of two dimensional multigrid}
According to Theorem 10.7.15 in \cite{Hackbusch1993}, the multigrid V-cycle converges for each of the two dimensional operators $\mathcal{L}_j^\hS$, $j=1,\dots,N$ if there exists a $C_A$ such that the smoothing property
\begin{equation}
  A^\hS_{\ell,j} \le W^\hS_{\ell,j}
  \label{eqn:SmoothingProperty2d}
\end{equation}
and the approximation property
\begin{equation}
  0 \le \left(A^\hS_{\ell+1,j}\right)^{-1} - P_\ell^\hS \left(A^\hS_{\ell,j}\right)^{-1} R_\ell^\hS \le C_A \left(W^\hS_{\ell+1,j}\right)^{-1}
  \label{eqn:ApproximationProperty2d}
\end{equation}
are satisfied on all levels $\ell=1,\dots,L$.

The smoothing property (\ref{eqn:SmoothingProperty2d}) is automatically satisfied for (sufficiently damped) point Jacobi and SOR smoothers (Remark 4.6.5 in \cite{Buckeridge2010PhD}). To see this, denote the matrix consisting only of the diagonal entries of $A_{\ell,j}^\hS$ by $D_{\ell,j}^\hS$ and use $W_{\ell,j}^\hS=\rho_{\operatorname{relax}}^{-1} D_{d,j}^\hS$, i.e. weighted point Jacobi relaxation. The relaxation parameter is chosen such that $\rho_{\operatorname{relax}} \le ||(D_{\ell,j}^\hS )^{-1}A_{\ell,j}^\hS||^{-1}$ where $||\cdot||$ is the spectral norm. Then (\ref{eqn:SmoothingProperty2d}) follows by definition from the equivalence \mbox{$-\Id\le X\le \Id$} $\Leftrightarrow$ $||X||\le 1$ applied to $X=(W_{\ell,j}^\hS)^{-1}A_{\ell,j}^\hS$.

A proof of the approximation property is significantly harder and we will not give it here (see Lemma 10.7.8 and Remark 10.7.13 in \cite{Hackbusch1993}). It depends on some minimal regularity assumptions on the profiles $\mat{\alpha}_\hS^\hS(\rhat)$ and $\alpha_r^\hS(\rhat)$. The constant $C_A$ may depend on the contrast, i.e. the maximum variation of the profiles.
We stress again that we use quasi-uniform grids for the horizontal discretisation (see the review in \cite{Staniforth2012} for a discussion of grids considered in meteorological application). In contrast to latitude-longitude grids, where the convergent grid lines near the pole introduce an additional horizontal anisotropy, the ratio between the smallest and largest grid spacing is bounded from below in the grids we consider. Hence the simple block-Jacobi and block-SOR smoothers which relax all degrees of freedom in one vertical column simultaneously will be efficient and no additional horizontal plane smoothing or selective semi-coarsening as described in \cite{Buckeridge2010} is required.

As the two dimensional equations are solved on the unit sphere, the operator $\mathcal{L}_j^\hS$ could become near-singular if $\lambda_j\rightarrow 0$. However, it is easy to see that this is not the case. As noted in Section \ref{sec:LinearEquation} we require the scaling $\omega\propto\Delta t\propto h_L$ to keep the Courant number fixed as the horizontal resolution increases. Therefore the second order term in (\ref{eqn:Cont2dOperator}) is of order $1$ and hence the relative importance of the two terms in (\ref{eqn:Cont2dOperator}) is independent of grid resolution.  It follows that all the eigenvalues $\lambda_j$ of (\ref{eqn:EigenvalueEquation}) are of order 1. It is a reasonable assumption that the profiles $\mat{\alpha}_\hS^\hS(\rhat)$, $\alpha_r^\hS(\rhat)$ are ``well-behaved'' in the sense that they are dominated by large scale variations due to global weather systems, small scale phenomena such as strong local variations carry substantially less energy. In this case we expect the spectrum of $\mathcal{L}_j^\hS$ to be bounded from above and below by two constants which are independent of $h_L$.

\paragraph{Convergence of three dimensional multigrid}
As argued above, the three dimensional problem can be decoupled into a set of $n_r$ two-dimensional problems. Due to the particular form of the smoother and of the prolongation/restriction matrices, it is in fact easy to verify that the smoothing property and the approximation property
\begin{xalignat*}{2}
  A^\otimes_{\ell} &\le W^\otimes_{\ell} \quad \text{and} &
  0 \le \left(A^\otimes_{\ell+1}\right)^{-1} - P_\ell^h \left(A^\otimes_{\ell}\right)^{-1} R_\ell^h &\le C_A \left(W^\otimes_{\ell+1}\right)^{-1},
\end{xalignat*}
for the tensor product multigrid algorithm for the original 3D problem on $\Omega$ follow directly from the respective properties (\ref{eqn:SmoothingProperty2d}) and (\ref{eqn:ApproximationProperty2d}) for the 2D problems on $\hS$, for all $j=1,\dots,N$.
\begin{theorem}
Let us assume that (\ref{eqn:SmoothingProperty2d}) and (\ref{eqn:ApproximationProperty2d}) are satisfied, for all $j=1,\dots,N$, and
let $\Mfac$ denote the iteration matrix for one step of the tensor product multigrid V-cycle defined above, i.e.
\begin{equation*}
  \mathbf{u} - \mathbf{u}_\otimes^* \ \mapsto \ \Mfac\left(\mathbf{u}-\mathbf{u}_\otimes^*\right)
\end{equation*}
where $\mathbf{u}_\otimes^*$ is the exact solution of the equation $\Afac \mathbf{u}^*_\otimes = \mathbf{f}$. Then the convergence rate
\begin{equation*}
 \rho_A^\otimes = \facNorm{\Mfac} \le \frac{C_A}{C_A+2(\nu_{\text{pre}}+\nu_{\text{post}})} \ < \ 1
\end{equation*}
independent of $h_L$, where $\facNorm{\cdot}$ is the energy norm induced by $\Afac$.
\end{theorem}

This is the main result given and proved for the two dimensional case in \cite[Theorem 2]{BoermHiptmair1999}. As we have seen above, the proof extends directly also to three dimensions and to our pressure correction problem here. In that case the assumptions of the theorem are satisfied as discussed above.

\subsubsection{Non-factorising case}
\label{sec:ProofNonFactorising}
We now assume that the matrix $A$ can be written as the sum of a perfectly factorising symmetric positive definite matrix $\Afac$ and a small correction $\deltaA$, namely $A = \Afac+\deltaA$. We quantify the deviation from perfect factorisation by $\Delta = \facNorm{\left(\Afac\right)^{-1}\deltaA}$ and assume that $\Delta < 1$. We also assume that the theory in Section \ref{sec:ProofFactorising} applies and the multigrid iteration for the factorising operator $\Afac$ converges, i.e. the error is reduced by a factor $\rho^\otimes_A <1$ in every multigrid V-cycle. The Richardson iteration for the full operator $A$ preconditioned with $\mu$ multigrid V-cycle cycles for $\Afac$ can then formally be written as
\begin{equation*}
  \mathbf{u}^{(k+1)} = \mathbf{u}^{(k)} + \left[\Id - \left(\Mfac\right)^\mu\right]\left(\Afac\right)^{-1}\left(\mathbf{f}-A\mathbf{u}^{(k)}\right).
\end{equation*}
Then at every step the error $\mathbf{u}^{(k)}-\mathbf{u}^*$ to the exact solution $\mathbf{u}^* := A^{-1} \mathbf{f}$ is reduced by a factor
\begin{align}
  \rho_A & \ = \ \facNorm{\Id - \left[\Id-\left(\Mfac\right)^\mu\right]\left(\Afac\right)^{-1}A} \nonumber\\
& \ \le \ \facNorm{\left(\Afac\right)^{-1}A-\Id}+\facNorm{\Mfac}^\mu\facNorm{\left(\Afac\right)^{-1}A}
\ \le \ \Delta + (1+ \Delta) (\rho^\otimes_A)^\mu.
  \label{eqn:RichardsonConvergenceRate}
\end{align}
Thus, for an arbitrary $\Delta < 1$ the convergence rate $\rho_A$ is less than 1, provided the number of V-cycles $\mu > \log_{\rho^\otimes_A} \left( (1-\Delta)/(1+\Delta) \right)$. On the other hand, if we only apply one V-cycle ($\mu = 1$), then a convergence rate $\rho_A < 1$ can still be
guaranteed provided $\Delta < (1-\rho^\otimes_A)/(1+\rho^\otimes_A)$.
Similar results can also be proved for the convergence of Krylov solvers, such as BiCGStab, preconditioned with $\mu$ multigrid V-cycle cycles for $\Afac$.
\section{Discretisation and Implementation}\label{sec:Implementation}
In practise, and as we demonstrate in the following, the tensor product preconditioners will be efficient for a wider range of problems not covered by the formal theory. We now describe the discretisation and DUNE implementation of the solvers we used in our numerical experiments.
\subsection{Grid structure and discretisation}
For simplicity we use a simple finite volume discretisation for all numerical experiments in this work. More complex schemes such as mimetic mixed finite elements are also currently under consideration for the development of dynamical cores \cite{Cotter2012,Cotter2013} and might require the solution of the equation in a different pressure space, such as higher order DG space. However, the basic ideas described in this work can still be applied.

Grids used in meteorological applications (and also in many ocean models \cite{Marshall1997,Fringer2006}) usually have a tensor-product structure. They consist of a semi-structured two dimensional horizontal grid on the surface of the sphere and a one-dimensional vertical grid which is often graded to achieve higher resolution near the surface. In particular each three dimensional grid cell $E=(T,k)$ can be uniquely identified by the corresponding horizontal cell $T$ and a vertical index $k\in 1,\dots,n_r$.
This tensor-product structure in itself has important implications for the performance of any implementation: while it might be necessary to use indirect indexing for the horizontal grid, the vertical grid can always be addressed directly. As typically the number of vertical levels is large with $n_r\gtrsim100$, the cost of indirect addressing in the horizontal direction can be ``hidden'' \cite{MacDonald2011}, a phenomenon which we have confirmed numerically for our solvers in Section \ref{sec:IndirectAddressing}.
Furthermore fields can be stored such that the levels in each column are stored consecutively in memory, which leads to efficient cache utilisation (however, as discussed in \cite{Mueller2013a} a different memory layout has to be used on GPU architectures where the vertically-consecutive storage would prevent coalesced memory access in the tridiagonal solve).
To be able to use the geometric multigrid solvers described in this work, we also assume that the horizontal grid has a natural hierarchy; this is true for the icosahedral grids which are used in our numerical tests where each triangular coarse grid cell consist of four smaller triangles on the next-finer multigrid level. In contrast to a simple longitude-latitude grid, these semi-structured grids have no pole-problem, i.e. the ratio between the size of the largest and smallest grid spacing is bounded. This implies that there is no additional horizontal anisotropy which would further complicate the construction of a solver (however, as has been shown in \cite{Buckeridge2010,Buckeridge2011}, the tensor-product multigrid approach can still be applied for longitude-latitude grids if the horizontal coarsening strategy is adapted appropriately).

In the finite volume discretisation any continuous field $u(r,\rhat)$ is approximated by its average value in a grid cell. In particular, for each \textit{horizontal} grid cell $T$ we store one vector $\mathbf{u}_T$ of length $n_r$ representing the field in the vertical column. In this cell the discrete equation (\ref{eqn:AlgebraicEquation}) for the $n_r$-vector $\mathbf{u}_T$ can be written as
\begin{equation}
  (\mat{A}\mathbf{u})_T = \mat{A}_T \mathbf{u}_T + \sum_{T'\in\nbT}\mat{A}_{TT'} \mathbf{u}_{T'} = \mathbf{f}_T,
\label{eqn:ColumnEquation}
\end{equation}
where the sum runs over all horizontal neighbours $T'\in\mathcal{N}(T)$ of $T$. In this expression $\mat{A}_T$ and $\mat{A}_{TT'}$ are $n_r\times n_r$ tridiagonal- and diagonal matrices of the form
\begin{xalignat}{2}
  A_T &= \operatorname{tridiag}(\mathbf{a}_T,\mathbf{b}_T,\mathbf{c}_T), &
  A_{TT'} &= \operatorname{diag}(\mathbf{d}_{TT'}).
  \label{eqn:TridiagonalSystem}
\end{xalignat}
Both matrices can be reconstructed on-the-fly from a number of scalar quantities, which are obtained from a discrete approximation of the profiles in (\ref{eqn:profiles}) and geometric factors. This reduces the amount of main memory access, in particular if the factorising profiles in the \TPMGfac\ preconditioner are used. For each horizontal cell $T$ the explicit expressions of the diagonals $\mathbf{a}_T$, $\mathbf{d}_{TT'}$ and upper- and lower- subdiagonals $\mathbf{b}_T$, $\mathbf{c}_T$ depend on whether the profiles can be factorised or not and are given explicitly in the next section. A block-SOR iteration with overrelaxation factor $\rho_{\operatorname{relax}}$ can then be written as
\begin{equation}
 \mathbf{u}_T \mapsfrom \mathbf{u}_T + \rho_{\operatorname{relax}}
  \left(\mat{A}_T\right)^{-1}\left(\mathbf{b}_T-\left(\mat{A}\mathbf{u}\right)_T\right)
  \label{eqn:BlockJacobi}
\end{equation}
and requires a tridiagonal solve in each vertical column to apply the inverse of the matrix $\mat{A}_T$ to the residual. This can be implemented using the Thomas algorithm \cite{Press2007}.
\subsection{Matrix-free DUNE Implementation}
All code was implemented using the DUNE library \cite{Bastian2008a,Bastian2008b}, which provides a set of C++ classes for solving PDEs using grid based methods. In particular it provides interfaces to (parallel) grid implementations such as ALUGrid \cite{Schupp1999,Dedner2004,Burri2005} and UGGrid \cite{Bastian1997}. The implementation of the grids is separated from data which is attached to the grid by the user via mapper functions between different grid entities (cells, edges, vertices) and the local data arrays. In our case we used the DUNE-grid module to implement a two dimensional host grid and then attached a whole column of length $n_r$ to each horizontal grid cell $T$. We represent the matrix as follows: in the non-factorising case (\TPMGfull), we store a vector $\hat{\vec{\beta}}_T$ of length $n_r$ at each horizontal cell to represent the zero order term, two vectors $(\hat{\vec{\alpha}}_r)_T$ and $(\hat{\vec{\xi}}_r)_T$ of length $n_r+1$ to represent the vertical diffusion and advection terms, and one vector  $(\hat{\vec{\alpha}}_\hS)_{TT'}$ of length $n_r$ at each horizontal edge $TT'$. The explicit form of these vectors is obtained by a standard finite volume discretisation of the problem%
\ifpreprint 
 and is written down in equation (\ref{eqn:hatFunctionsNonFac}) in the appendix
\fi 
. The vectors $\mathbf{a}_T$, $\mathbf{b}_T$, $\mathbf{c}_T$ and $\mathbf{d}_{TT'}$ in (\ref{eqn:TridiagonalSystem}) are
\begin{equation}
  \begin{aligned}
    d_{TT',k} &= -(\hat{\alpha}_{\hS})_{TT',k},&
    d_{T,k} &= \sum_{T'\in\nbT} \tilde{d}_{TT',k},\\
    b_{T,k} &= -(\hat{\alpha}_r)_{T,k+1}-(\hat{\xi}_r)_{T,k+1},&
    c_{T,k} &= -(\hat{\alpha}_r)_{T,k+1}+(\hat{\xi}_r)_{T,k+1},\\
    a_{T,k} &= \hat{\beta}_{T,k} - (b_{T,k}+c_{T,k} + d_{T,k}).
  \label{eqn:TridiagonalNonFactorising}
  \end{aligned}
\end{equation}
In the factorising case (\TPMGfac) it is only necessary to store \textit{scalars} $\hat{\beta}^\hS_T$, $(\hat{\alpha}_r^\hS)_T$, $(\hat{\xi}_r^\hS)_T$ and $(\hat{\alpha}_\hS^\hS)_{TT'}$  on the horizontal cells and edges. In addition to this, four vectors of length $n_r$ and $n_r+1$ ($\hat{\vec{\beta}}^r$, $\hat{\vec{\alpha}}_\hS^r$, $\hat{\vec{\alpha}}_r^r$ and $\hat{\vec{\xi}}_r^r$) which arise from the vertical discretisation need to be stored once for the entire grid.
\ifpreprint 
The explicit form of these quantities is given in (\ref{eqn:hatFunctionsFac}) in the appendix.
\fi 
Similarly to (\ref{eqn:TridiagonalNonFactorising}) the matrix entries in (\ref{eqn:TridiagonalSystem}) can be calculated on the fly as
\begin{equation}
  \begin{aligned}
    d_{TT',k} &= -(\hat{\alpha}^r_{\hS})_{k}(\hat{\alpha}^\hS_{\hS})_{TT'},&d_{T,k} &= \sum_{T'\in\nbT} \tilde{d}_{TT',k}\\
    b_{T,k} &= -(\hat{\alpha}^r_r)_{k+1}(\hat{\alpha}^\hS_r)_{T}-(\hat{\xi}^r_r)_{k+1}(\hat{\xi}^\hS_r)_{T},&
    c_{T,k} &= -(\hat{\alpha}^r_r)_{k+1}(\hat{\alpha}^\hS_r)_{T}+(\hat{\xi}^r_r)_{k+1}(\hat{\xi}^\hS_r)_{T},\\
    a_{T,k} &= \hat{\beta}^r_{k}\hat{\beta}^\hS_{T} - (b_{T,k}+c_{T,k} + d_{T,k}).
  \end{aligned}
  \label{eqn:TridiagonalFactorising}
\end{equation}
The scalars $\hat{\beta}^\hS_T$, $(\hat{\alpha}_r^\hS)_T$, $(\hat{\xi}_r^\hS)_T$ and $(\hat{\alpha}_\hS^\hS)_{TT'}$ only need to be read once per vertical column and the associated cost can be hidden together with the cost of indirect addressing on the horizontal grid for large enough $n_r$. Moreover, the vectors $\hat{\vec{\beta}}^r$, $\hat{\vec{\alpha}}_\hS^r$, $\hat{\vec{\alpha}}_r^r$ and $\hat{\vec{\xi}}_r^r$ require only a small amount of memory and can be cached. In summary, the cost of memory access for the matrix is likely to be significantly smaller than the cost of accessing field vectors such as $\mathbf{u}_T$ and $\mathbf{b}_T$ when solving the tridiagonal system in (\ref{eqn:BlockJacobi}) or in the matrix vector product.

The DUNE-grid interface provides iterators over the horizontal grid cells and over the neighbours of each cell. To implement for example the sparse matrix vector product (SpMV) in (\ref{eqn:ColumnEquation}) we iterate over all horizontal grid cells $T$, and then in each cell we loop over the edges $TT'$ for all neighbours $T'$ to read the profiles stored on the cells and edges from memory and construct the matrices $\mat{A}_T$ and $\mat{A}_{TT'}$. These are then applied to the local vectors $\mathbf{u}_T$ and $\mathbf{u}_{T'}$ to evaluate $(\mat{A}\mathbf{u})_T$, which requires inner loops over the vertical levels.
Of all grids that are currently available through the DUNE interface we found that only ALUGrid can be used to represent a two-dimensional sphere embedded in three dimensional space. Unfortunately the scalability of ALUGrid is very limited because in a parallel implementation the entire grid is stored on each processor. Alternatively we used a three dimensional UGGrid implementation for a thin spherical shell consisting of one vertical layer to represent the unit sphere. Based on the coarsest grid, finer multigrid levels can be constructed by refinement in the horizontal direction only. Any geometric quantities in this thin three dimensional grid can then be related to the corresponding values on the two dimensional grid by simple scaling factors. We implemented both a gnomonic cubed sphere grid \cite{Sadourny1972} and an icosahedral grid, for which the grid points are projected onto the sphere, and all numerical results reported in this work were obtained with the icosahedral grid.

As is typical in atmospheric applications, parallel domain decomposition is in the horizontal direction only. As the DUNE host grids that we used are already inherently parallel, parallelisation of the code was straightforward by calling the relevant halo exchange routines when necessary. Load balancing was achieved by choosing the problem size such that the number of cells on the coarses level is identical to the number of processors and each processor ``owns'' one coarse grid cell and the corresponding child cells.
While at first sight this might cause a problem for large core counts because the coarsest level still has a relatively large number of degrees of freedom and the multigrid hierarchy is very shallow, it turns out that the zero order term in the Helmholtz equation (\ref{eqn:HelmholtzVector}) averts potential problems. This is because relative to the zero order term the importance of the horizontal diffusion term decreases with a factor of four on each coarse level, and so after a small number of coarsening steps the problem is well conditioned and can be solved by a very small number of smoothing iterations. An alternative and more physical explanation is that any interactions in the continuous PDE in (\ref{eqn:HelmholtzVector}) are exponentially damped with an intrinsic length scale $\omega$ and hence it is not necessary to coarsen the grid beyond this scale. This has been confirmed numerically for a simplified test problem in \cite{Mueller2013}, where it has been shown that as little as four multigrid levels still give very good convergence for typical grid spacings and time step sizes. In the parallel scaling tests in this work we typically used 6 or 7 multigrid levels and one iteration of the smoother to solve the coarse grid problem.
\section{Numerical results}\label{sec:NumericalResults}
In the following we study the performance of the two tensor-product preconditioners \TPMGfull\ and \TPMGfac\ described in Section \ref{sec:TPMGs} applied to two test cases in atmospheric flow simulation. We confirm the optimality and robustness of \TPMGfull\ even for non-factorising profiles, compare the performance of the two variants and study their parallel scalability. All runs (including the sequential tests) were carried out on the phase 3 configuration of the HECToR supercomputer, which consists 2816 compute nodes with two 16-core AMD Opteron 2.3GHz Interlagos processors each. The entire cluster contains 90,112 cores in total. The code was compiled with version 4.6.3 of the gnu C compiler.

Unless stated otherwise we always used 6 multigrid levels with two vertical line-SOR pre- and post- smoothing steps on each level ($\nu^{\text{pre}}=\nu^{\text{post}}=2$); the overrelaxation factor in the smoother was set to $\rho_{\text{relax}}=1$. One smoother iteration is used to solve the coarse grid problem. We use linear interpolation to prolongate the solution to the next-finer grid ($2h\mapsto h$). The right hand side, which in each cell represent a cell integral of a field, is restricted to the next-coarser level ($h\mapsto 2h$) by summing the fine grid values of all four fine grid cells comprising the coarse grid cell. Recall that these integrid-operations only require interpolation and summation in the horizontal direction. The tolerance in the iterative solver was set to $10^{-5}$, i.e. we iterate until the residual has been reduced by at least five orders of magnitude. The number of vertical levels was set to $n_r=128$, which is typical for current meteorological applications. We note, however, that all runtimes should be directly proportional to $n_r$ (and this is confirmed in the following section).
\subsection{Overhead from indirect addressing}
\label{sec:IndirectAddressing}
\begin{figure}
\begin{center}
  \includegraphics[width=0.5\linewidth]{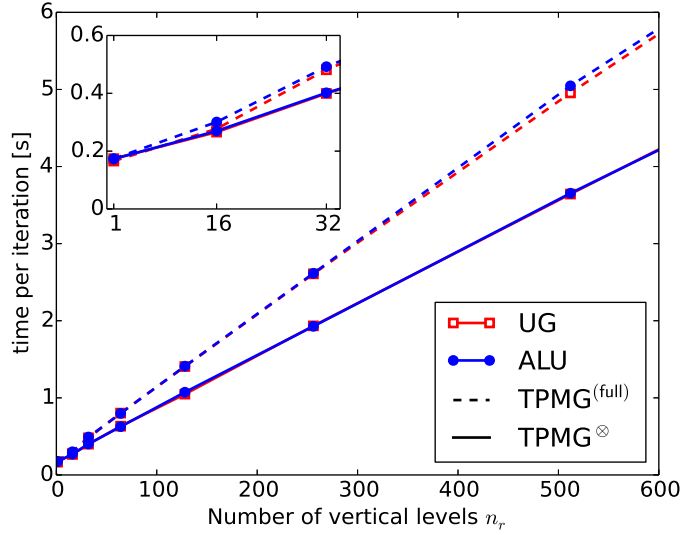}
  \caption{Time per iteration for different numbers $n_r$ of vertical levels for two grid implementations (UGGrid in red, open squares and ALUGrid in blue, filled circles) on an icosahedral grid; results for shown both for the \TPMGfull\ (dashed lines) and \TPMGfac\ (solid lines) preconditioner.}
  \label{fig:nzDependency}
\end{center}
\end{figure}

While data in one vertical column is stored consecutively in memory and can be addressed directly, in general indirect addressing has to be used in the horizontal directions. However, as the horizontal lookup is only required once per column, the relative penalty for this will be very small provided $n_r$ is large enough. As discussed in \cite{MacDonald2011}, in this case the overhead from indirect addressing can be ``hidden'' behind work in the vertical direction. To verify this we ran our solver with two different DUNE grid implementations and measured the time per iteration for different numbers of vertical levels. We expect this time to depend on $n_r$ as follows
\begin{equation}
 t_{\operatorname{iter}}(\operatorname{grid},n_r) = \left(C_0+C_{\operatorname{grid}}\right) + q\cdot n_r
\label{eqn:nzDependency}
\end{equation}
where $C_{\operatorname{grid}}$ is the overhead of indirect addressing and depends on the grid implementation. The constant $C_0$ encapsulates any other work which is only done once per column and both $C_0$ and the slope $q$ are independent of the horizontal grid. Figure \ref{fig:nzDependency} shows the results for the ALUGrid and UGGrid implementation and confirms the linear dependency in (\ref{eqn:nzDependency}).
As can be seen from this plot, for both preconditioners \TPMGfull\ and \TPMGfac\ the overhead from indirect addressing $C_{\operatorname{grid}}$ and the additional overhead $C_0$ together are at the order of less than $20\%$ as soon as $n_r\gtrsim 100$. Incidentally both DUNE grid implementations that we tested are equally efficient. We stress that in both grids data in adjacent vertical columns is not necessarily stored consecutively in memory. Not surprisingly, the slope $q$ is larger for the more expensive $\TPMGfull$ preconditioner. The results in this section also confirm that performance tests carried out on a directly addressed horizontal grid, such as the results in \cite{Mueller2013}, can be generalised to indirectly addressed grids.
\subsection{Test Case I: Balanced zonal flow}
\label{sec:BalancedFlow}
\begin{figure}
  \begin{minipage}{0.425\linewidth}
  \includegraphics[width=0.95\linewidth]{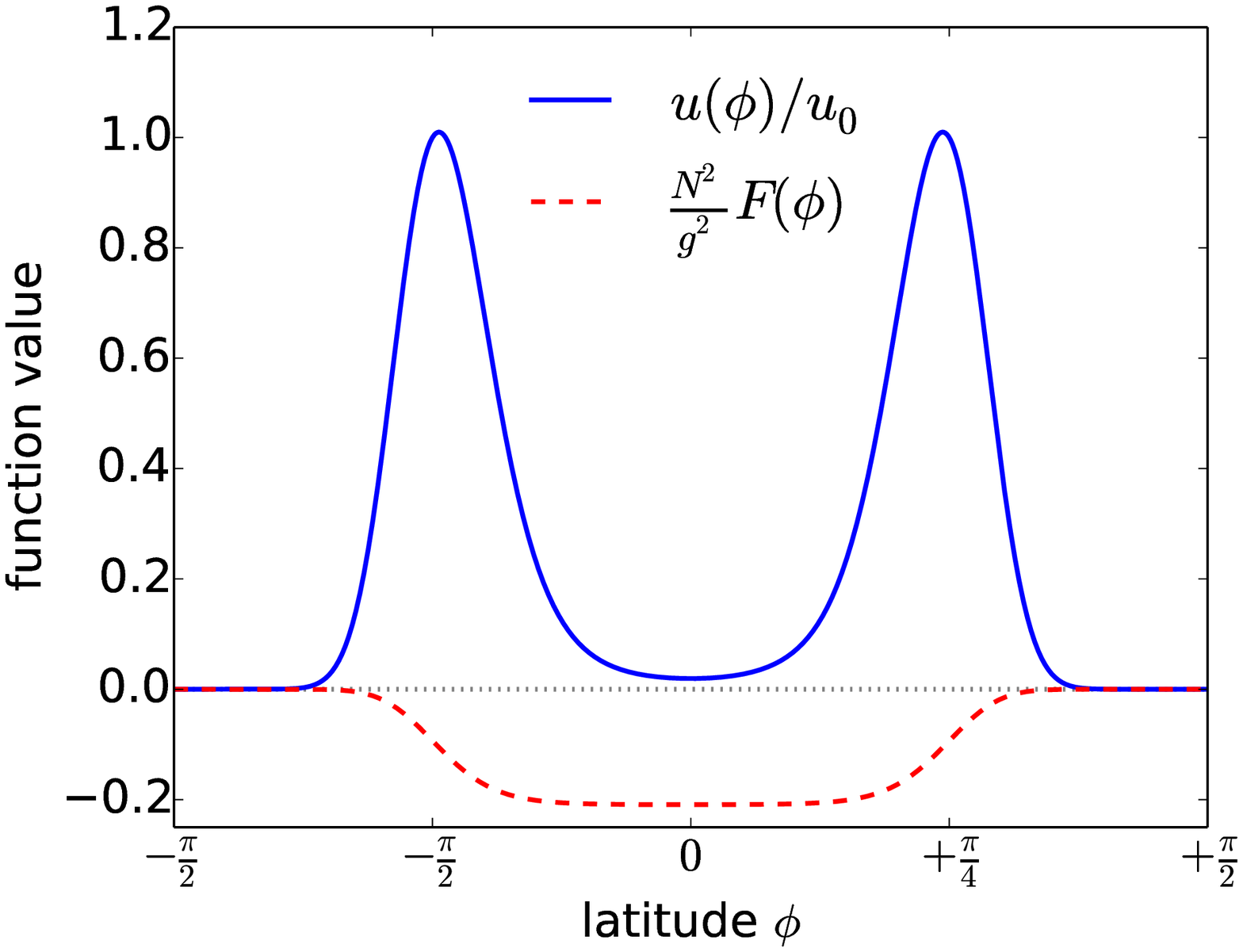}
  \end{minipage}
  \hfill
  \begin{minipage}{0.55\linewidth}
    \includegraphics[width=0.95\linewidth]{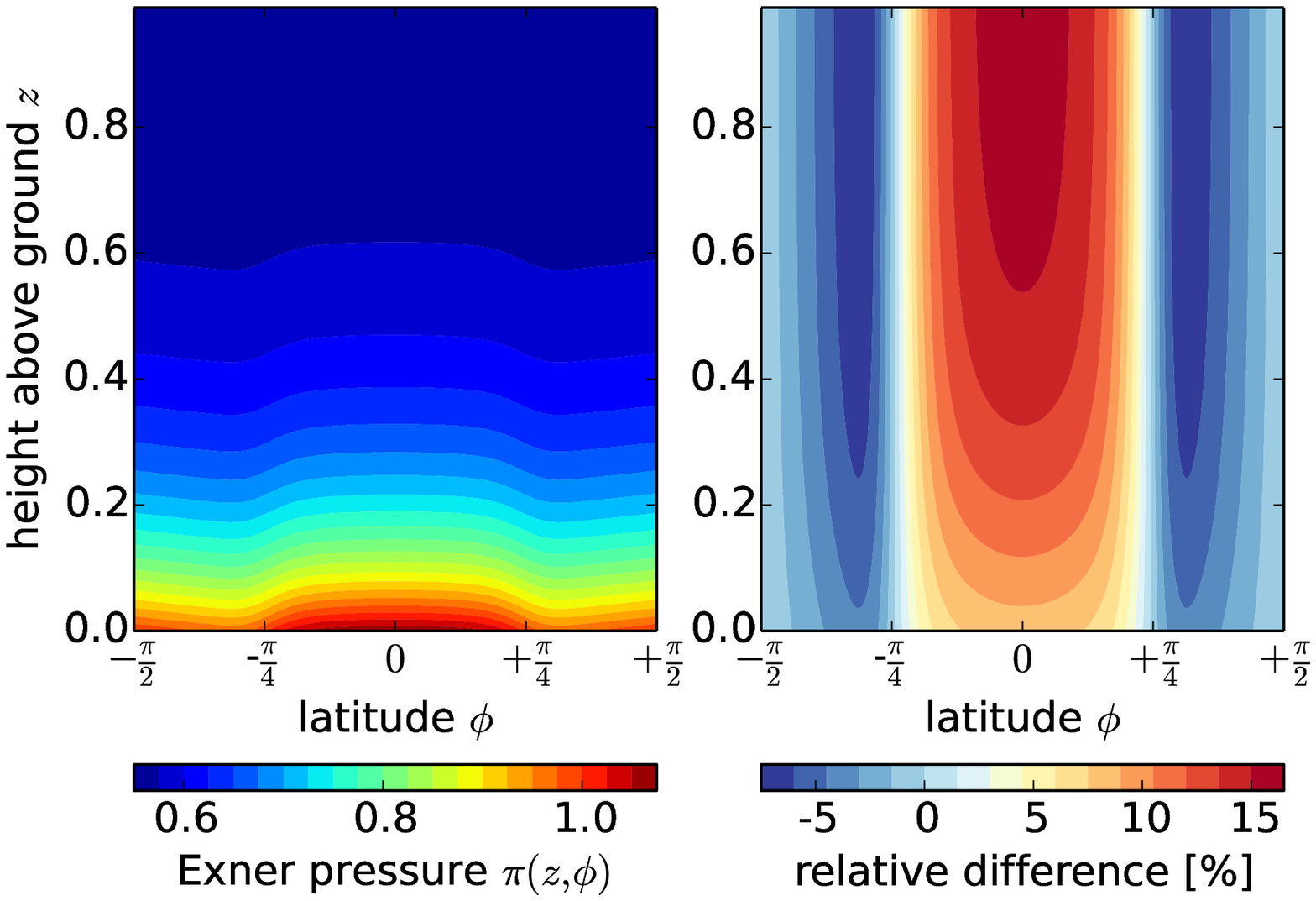}
  \end{minipage}
  \caption{Left: Velocity field $u(\phi)$ and jet function $F$ defined in eqns. (\ref{eqn:JetFunction}) and (\ref{eqn:VelocityField}) for $N=0.028s^{-1}$ ($\epsilon=1.23$). Right: Exner pressure $\pi$ and relative difference $\frac{\pi^\otimes-\pi}{\pi}$ in the $(\phi,z)$-plane for the same value of $N$. The height above ground is measured in units of the depth of the atmosphere.}
    \label{fig:ExnerPressure}
\end{figure}

We first test our solver with the profiles from a simplified meteorological test problem which corresponds to a balanced atmosphere with constant buoyancy frequency and zonal flow with one jet in each hemisphere. The advantage of this test case is that the deviation of the atmospheric profiles from a perfect factorisation can be controlled by varying a single parameter. In \cite{Davies2005} it is shown that under the assumption that the velocity field points in the longitudinal direction and the buoyancy frequency $N$ defined in (\ref{eqn:BuoyancyFrequency}) is constant, a solution of the Euler equations is given by
\begin{equation}
  \begin{aligned}
  \pi(\rhat,r) &= \frac{\epsilon+\eHoriz\eVert}{1+\epsilon},&
  \theta(\rhat,r) &= T_0 \left(\eHoriz\eVert\right)^{-1}, \\
  \rho(\rhat,r) &= \frac{p_0}{R_dT_0} \left[\pi(\rhat,r)\right]^{\gamma}\eHoriz\eVert,&
  u(\rhat,r) &= u_\hS(\phi)
  \end{aligned}
  \label{eqn:BalancedFlow}
\end{equation}
where the functions $\eHoriz$ and $\eVert$ are defined as
\begin{xalignat*}{2}
  \eHoriz &= \exp\left[-\frac{N^2}{g^2}F(\phi)\right], &
  \eVert &= \exp\left[-\frac{N^2\Rearth}{g}(r-1)\right].
\end{xalignat*}
In the horizontal direction the profiles only vary in the latitudinal direction  $\phi\in[-\pi/2,\pi/2]$. The parameter $\epsilon$ is related to the buoyancy frequency by $\epsilon \equiv \left(\frac{N}{N^*}\right)^2-1$ with $N^* \equiv \frac{\sqrt{c_pT_0}}{g} = \frac{c_h}{g}$.
The function $F(\phi)$ is related to the velocity field $u_\hS(\phi)$ as
\begin{equation}
  \frac{dF(\phi)}{d\phi} = 2\Rearth\Omegaearth u_\hS(\phi) \sin\phi + u_\hS(\phi)^2 \tan\phi
  \label{eqn:JetFunction}
\end{equation}
with angular velocity $\Omegaearth=2\pi/(24\cdot 3600)s^{-1}$.
For our numerical experiments we choose the velocity such that it corresponds to two jets with peak velocity $u_0=100ms^{-1}$ in the mid latitudes ($\phi_M=\pi/4$, $\sigma=0.1$):
\begin{equation}
  u_\hS(\phi) = u_0 \frac{\cos\phi}{\cos\phi_M} \exp\left[-\frac{(\cos\phi-\cos\phi_M)^2}{2\sigma^2}\right]
  \label{eqn:VelocityField}
\end{equation}
as plotted together with the corresponding $F(\phi)$ in Figure \ref{fig:ExnerPressure}.
If we fix the reference pressure and temperature to physically realistic values $p_0=10,000 Pa$ and $T_0=273K$, the only free parameter in (\ref{eqn:BalancedFlow}) is the buoyancy frequency. In particular if $N$ is identical to $N^*$, i.e. $\epsilon=0$, the first term in the expression for the Exner pressure in (\ref{eqn:BalancedFlow}) vanishes and all profiles factorise exactly.

In the following we present numerical results for a range of buoyancy frequencies between $N=N^*=0.01873s^{-1}$ and $N=0.028s^{-1}$.
As a preconditioner we use both a multigrid algorithm with the full model operator and the tensor-product multigrid algorithm with an approximate factorisation of the Exner pressure
\begin{equation}
  \pi^{\otimes}(\rhat,r) = \pi^\hS(\rhat)\pi^r(r)
  \equiv \frac{\epsilon+\eVert}{1+\epsilon}\cdot\eHoriz
  \label{eqn:FacApproxBalancedFlow}
\end{equation}
which reduces to the expression in (\ref{eqn:BalancedFlow}) for $\epsilon=0$.
Both the Exner pressure and the relative difference $\frac{\pi^\otimes-\pi}{\pi}$, which is an indicator of the quality of the factorisation, are plotted for $N=0.028s^{-1}$ in the $(z,\phi)$ plane in Figure \ref{fig:ExnerPressure}. As can be seen from this figure, the relative difference between the profiles can be larger than 15\%.

\begin{figure}
  \begin{minipage}{0.45\linewidth}
  \includegraphics[width=1.0\linewidth]{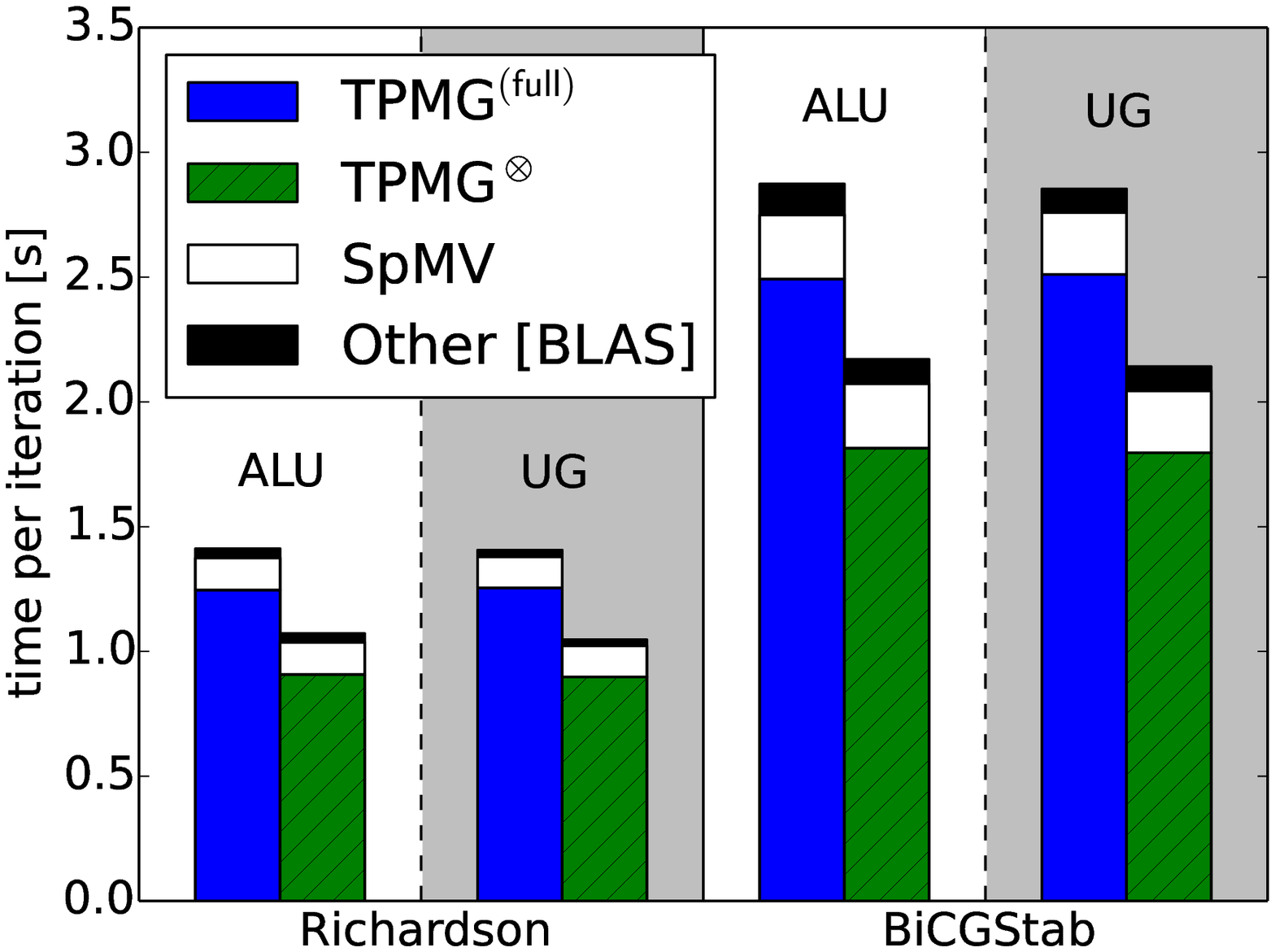}
  \end{minipage}
  \hfill
  \begin{minipage}{0.45\linewidth}
    \includegraphics[width=1.0\linewidth]{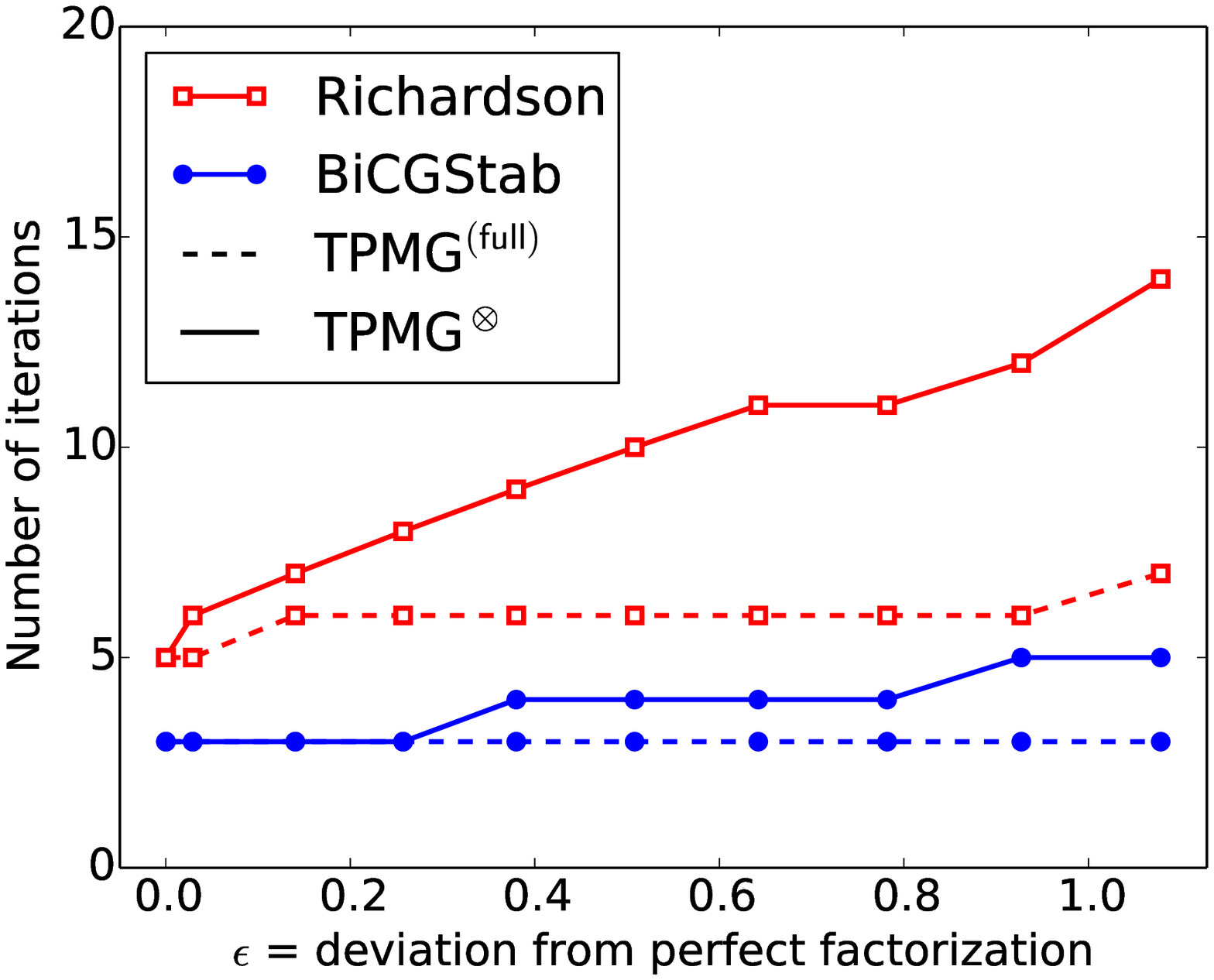}
  \end{minipage}
  \caption{Breakdown of the time per iteration for two different iterative solvers and grid implementations (left) and number of iterations (right) to reduce the relative residual by at least five orders of magnitude for the idealised balanced flow testcase. The multigrid preconditioner is used with both the full, non-factorising profiles (\TPMGfull, blue columns and dashed curves) and the approximate factorisation (\TPMGfac, hatched green columns and solid curves) in (\ref{eqn:FacApproxBalancedFlow}). In all cases a problem with $n_r=128$ and $2.6\cdot10^{6}$ total degrees of freedom was solved sequentially on HECToR.}
  \label{fig:tIternIterBreakdown}
\end{figure}
The time per iteration is shown in Figure \ref{fig:tIternIterBreakdown} (left) for two grid implementations. Both a preconditioned Richardson iteration and BiCGStab are used with one multigrid V-cycle as a preconditioner. It is important to note that BiCGStab requires two applications of the preconditioner and two sparse matrix-vector products per iteration, while the Richardson iteration only requires one of each, and not surprisingly the figure demonstrates that most of the time is taken up by the multigrid preconditioner in all cases. The number of iterations for each of the combinations is plotted in Figure \ref{fig:tIternIterBreakdown} (right) for a range of $\epsilon$.

First of all we note the almost perfect robustness of the full preconditioner \TPMGfull\ for this test problem where the profiles strongly deviate from the factorising case, but the convergence of preconditioned Richardson iteration and preconditioned BiCGStab are essentially not affected. The practically observed convergence rate for the V-cycle (in the Richardson iteration) is around $\rho_A = 0.1$. This confirms the theoretical results in Sections \ref{sec:ProofFactorising} and \ref{sec:ProofNonFactorising}. BiCGStab converges in approximately half the number of iterations than Richardson, as expected.
In terms of time per iteration, the multigrid preconditioner with factorised profiles (\TPMGfac) can be up to $25\%$ faster than the algorithm with non-factorising profiles (\TPMGfull). However, this comes at the expense of an increase in the number of iterations for larger values of $\epsilon$
that can be seen in Figure \ref{fig:tIternIterBreakdown} (right).
While for the Richardson iteration the increase is almost threefold if \TPMGfac\ is used, this is much less dramatic for BiCGStab where \TPMGfac\ only requires twice as many iterations as \TPMGfull\ for the largest $\epsilon$.

\begin{figure}
  \begin{center}
    \includegraphics[width=0.45\linewidth]{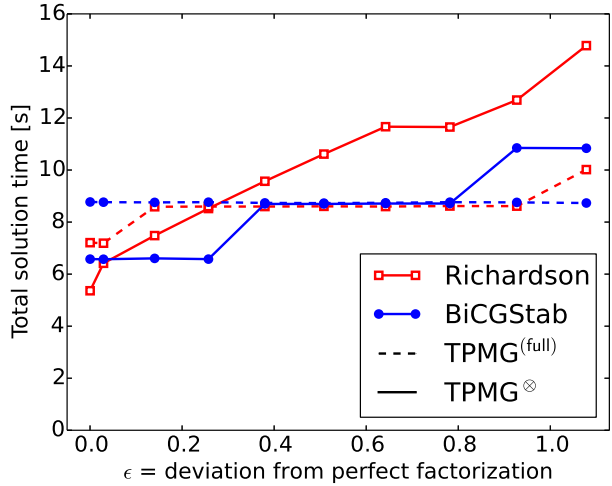}
    \caption{Total time required to reduce the relative residual by at least five orders of magnitude for the idealised balanced flow testcase. The multigrid preconditioner is used with both the full, non-factorising profiles (\TPMGfull, dashed curves) and the approximate factorisation (\TPMGfac, solid curves) in (\ref{eqn:FacApproxBalancedFlow}). In all cases a problem with $n_r=128$ and $2.6\cdot10^{6}$ total degrees of freedom was solved sequentially on one node of the HECToR supercomputer.}
    \label{fig:tSolveBalancedFlow}
  \end{center}
\end{figure}

Finally, the total solution time is shown in Figure \ref{fig:tSolveBalancedFlow}.
As expected, the total solution time for solvers with \TPMGfac\ preconditioner grows as $\epsilon$ increases. However, as the time per iteration is about 25\% smaller for this preconditioner, for small $\epsilon$ the total solution time is also reduced by a similar factor. The most robust solver appears to be BiCGStab, which gives the best overall performance for large $\epsilon$, even with the factorising preconditioner \TPMGfac.

\subsection{Test Case II: Aquaplanet}\label{sec:ResultsAquaplanet}
\begin{figure}
 \begin{minipage}{0.45\linewidth}
  \includegraphics[width=1.0\linewidth]{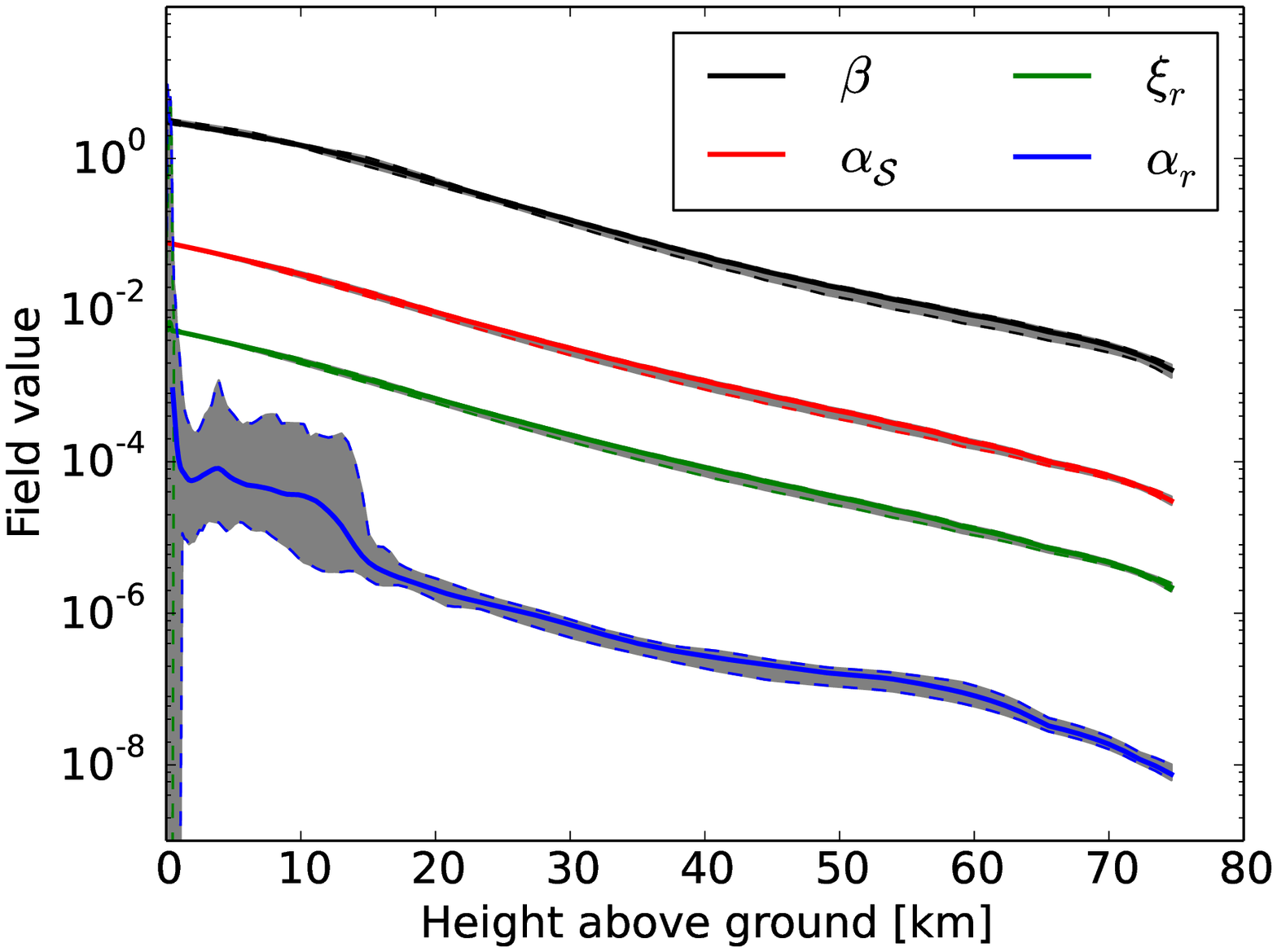}
  \end{minipage}
  \hfill
  \begin{minipage}{0.45\linewidth}
  \includegraphics[width=1.0\linewidth]{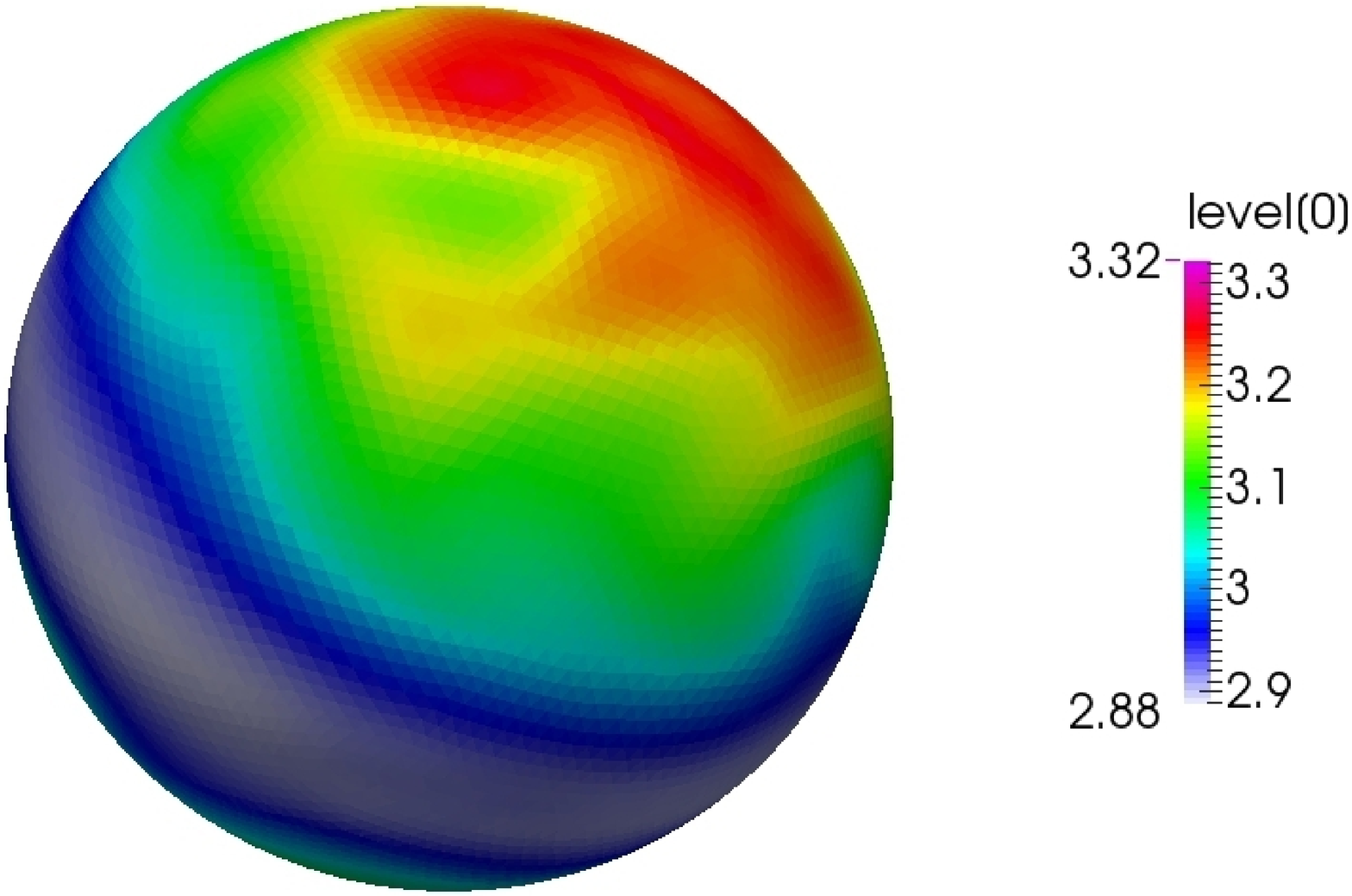}
  \end{minipage}
  \caption{Left: Dependency of the horizontal average of different profiles with height (solid lines, from top to bottom: $\beta$, $\alpha_\hS$, $\xi_r$ and $\alpha_r$). The horizontal variation is also represented by gray bands between the minimum and maximum value on each grid level (dashed curves). Right: Zero-order term $\beta=\gamma\rho/\pi$ on the lowest grid level. The horizontal variation in the field is at the order of $10\%$.}
  \label{fig:VerticalProfiles}
\end{figure}

While the runs in the previous section were carried out under idealised and not necessarily realistic conditions, we also tested our solver for profiles obtained from common meteorological test cases. We first obtained the profiles $\pi$, $\theta$ and $\rho$ from an aquaplanet run of the Met Office Unified Model\footnote{For technical reasons we used the wet density, such that the equation of state is not satisfied, but this should not have a significant impact on our conclusions.}. While these fields contain significantly more variation than the idealised profiles described in Section \ref{sec:BalancedFlow} and also describe phenomena such as convection near the ground and baroclinic instabilities, they are largely ``well behaved'' in the sense that most of them can be factorised approximately into a horizontal and a vertical variation. To quantify this further, we plot for each of the profiles the average, minimum and maximum over the horizontal grid on each vertical level in Figure \ref{fig:VerticalProfiles} (right).
For most profiles the horizontal variation is small and the average value decays exponentially with height; see for example Figure \ref{fig:VerticalProfiles} (left), which shows the profile $\beta$ on the lowest grid level. The only exception is $\alpha_{r}$ which shows significant horizontal variation in the lower atmosphere. This is mainly due to the fact that, as can be seen from the explicit expressions in (\ref{eqn:profiles}), this profile contains the buoyancy frequency and hence vertical derivatives of the potential temperature, which can vary significantly from column to column due to convection in the lower atmosphere.

We found that for these more typical profiles the factorising preconditioner \TPMGfac\ causes both solvers to diverge. An easy fix for this is to factorise all profiles except $\alpha_r$. We denote the resulting preconditioner with partial factorisation, where we keep the full non-factorising profile for $\alpha_r$, as \TPMGpartfac. As Table \ref{tab:tIterAquaplanet} demonstrates, this increases the time per iteration by just over $5\%$ relative to the fully factorising case (\TPMGfac), but it is still significantly smaller than in the non-factorised case (\TPMGfull).
\begin{table}
    \caption{Time per iteration and speedups relative to \TPMGfull\ for different solvers and preconditioners. In all cases a problem with $n_r=128$ and $2.6\cdot10^{6}$ total degrees of freedom was solved sequentially on one node of HECToR using the ALUGrid implementation.}
  \begin{center}
    \tabsize
    \label{tab:tIterAquaplanet}
    \begin{tabular}{lrrrrrr}
    \toprule
    Solver & \multicolumn{2}{c}{\TPMGfull} & \multicolumn{2}{c}{\TPMGpartfac} & \multicolumn{2}{c}{\TPMGfac} \\
    & $t_{\operatorname{iter}}$ & speedup
    & $t_{\operatorname{iter}}$ & speedup
    & $t_{\operatorname{iter}}$ & speedup\\
    \midrule
    Richardson & 1.43 & --- & 1.11 & $1.29\times$ & 1.05 & $1.36\times$ \\
    BiCGStab   & 2.88 & --- & 2.26 & $1.27\times$ & 2.14 & $1.35\times$ \\
    \bottomrule
    \end{tabular}
  \end{center}
\end{table}

The numbers of iterations and total solution times are shown in Table \ref{tab:ResultsAquaplanet}. We find that again the tensor product multigrid method converges extremely fast and robustly even though the profiles do not factorise, needing no more than 5 to 7 V-cycles to reduce the residual by 5 orders of magnitude. In this case the problem is solved fastest with the BiCGStab solver and the \TPMGpartfac\ preconditioner.  In total we find that, as in the idealised test case with small $\epsilon$, the (partially) factorised multigrid preconditioner can again lead to performance gains. As outlined in the introduction, these gains may be more significant on novel manycore architectures, such as GPUs, where the cost of memory references relative to one floating point operation is even larger.
\begin{table}
  \caption{Performance of different solvers for an aquaplanet run. A problem with $n_r=128$ and $2.6\cdot10^{6}$ total degrees of freedom was solved sequentially on HECToR using the ALUGrid implementation.}
  \label{tab:ResultsAquaplanet}
  \begin{center}
  \tabsize
  \begin{tabular}{lrrrr}
    \toprule
    & \multicolumn{2}{c}{\# iterations ($||r||/||r_0||$)} & \multicolumn{2}{c}{total time} \\
    Solver &
    \TPMGfull & \TPMGpartfac &
    \TPMGfull & \TPMGpartfac \\
    \midrule
    Richardson & 5 ($9.1\cdot10^{-6}$) & 7 ($4.8\cdot10^{-6}$) & 7.94 & 7.28 \\
    BiCGStab   & 3 ($5.1\cdot10^{-7}$) & 3 ($5.2\cdot10^{-6}$) & 8.81 & 6.94\\
    \bottomrule
    \end{tabular}
  \end{center}
\end{table}
\subsection{Parallel scaling tests}
\begin{figure}
 \begin{minipage}{0.45\linewidth}
  \includegraphics[width=1.0\linewidth]{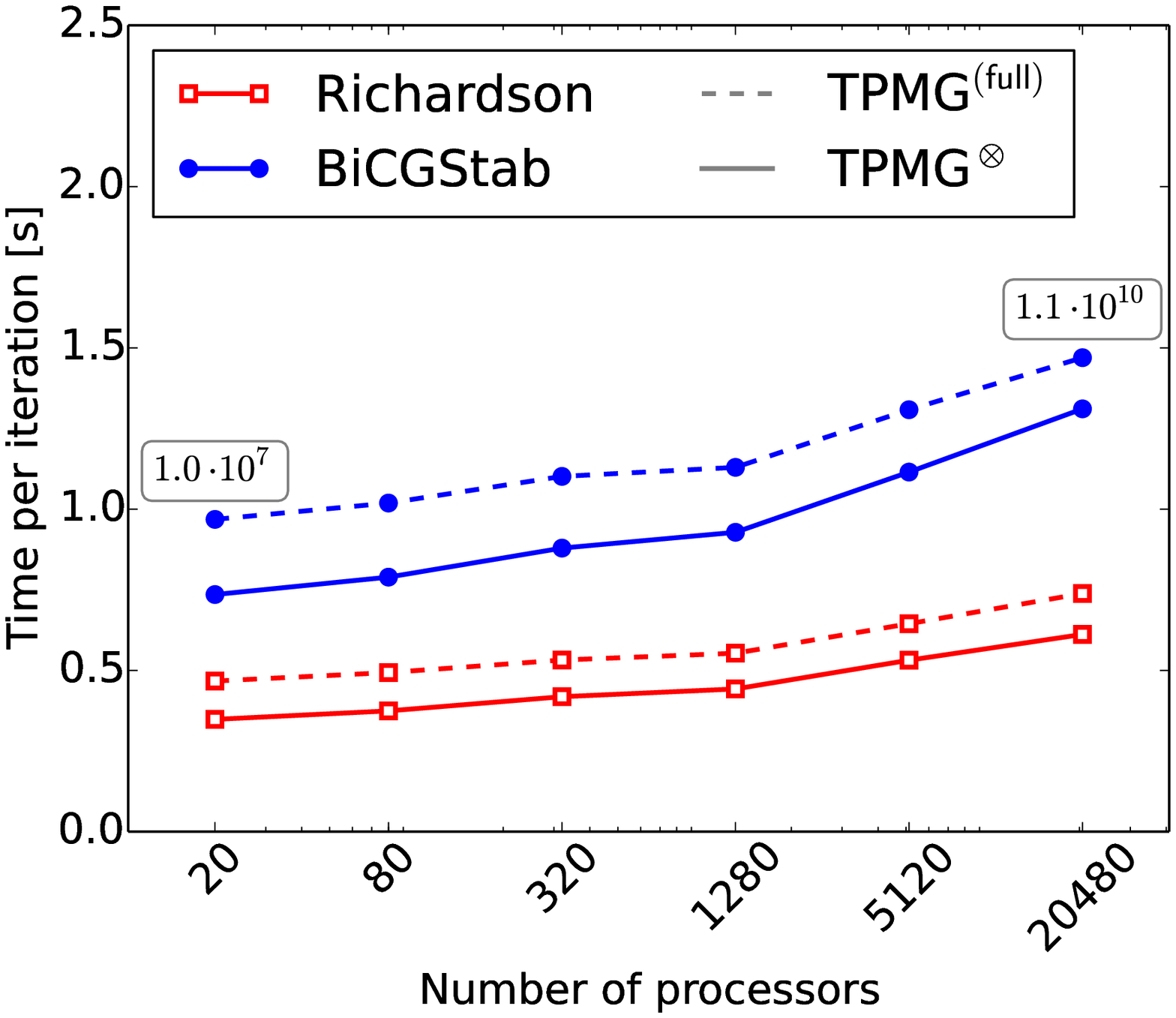}
 \end{minipage}
 \hfill
 \begin{minipage}{0.45\linewidth}
  \includegraphics[width=1.0\linewidth]{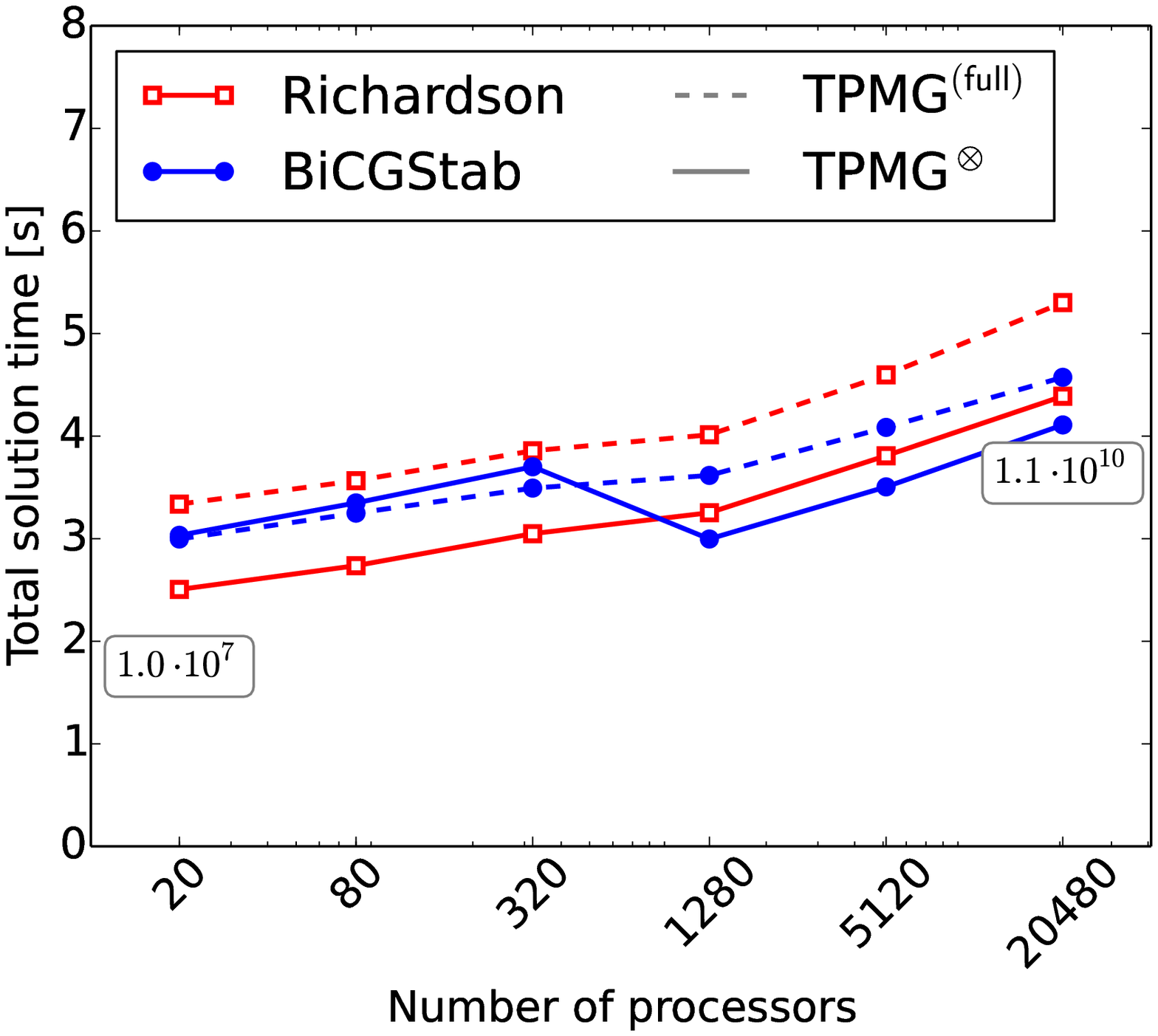}
 \end{minipage}
  \caption{Weak scaling of the time per iteration (left) and total solution time (right) on the HECToR supercomputer. The number of degrees of freedom varies from 10 million to 11 billion.}
  \label{fig:weakscaling}
\end{figure}

In addition to studying the sequential performance of the solvers, and in particular ensuring that they are algorithmically efficient, it is crucial to guarantee their parallel scalability on large computer clusters. For this we carried out scaling tests of our solvers for the balanced flow testcase described in Section \ref{sec:BalancedFlow} with $\epsilon=0.14$; in contrast to the previous runs we always used 7 multigrid levels so that on the coarsest level each processor stores one vertical column of data. In Figure \ref{fig:weakscaling} (left) the weak scaling of the time per iteration on the HECToR supercomputer is shown for up to 20,480 cores, the largest problem that was solved has just over $10^{10}$ degrees of freedom.
We find that the number of iterations does not increase with the core count, and even drops in some cases. The Richardson solver requires seven iterations to reduce the residual by five orders of magnitude for both preconditioners, whereas BiCGStab requires 4 (\TPMGfac) and 3 (\TPMGfull) iterations for the same residual reduction. Consequently the total solution time in Figure \ref{fig:weakscaling} (right) shows the same excellent weak scaling.
\section{Conclusion}\label{sec:Conclusions}
In this work we discussed several multigrid preconditioners for anisotropic problems in flow simulations in ``flat'' domains with high aspect ratio. The algorithms are based on the tensor-product multigrid approach proposed and analysed for two-dimensional problems with separable coefficients in \cite{BoermHiptmair1999}. We extended the method and its analysis to three dimensional problems and via a perturbation argument also to non-separable coefficients. We demonstrated the excellent performance of tensor-product multigrid for two model PDEs arising in semi-implicit semi-Lagrangian time stepping in atmospheric modelling. The numerical tests confirm the theoretically predicted optimality and effectivity of the method. The practically observed convergence rates are around $\rho_A = 0.1$. The tests also show that under certain conditions a preconditioner based on an approximate factorisation of the atmospheric profiles can reduce the total solution time. We found this to be the case both for an idealised flow scenario and for a more realistic aquaplanet test case. We also demonstrated the excellent weak parallel scaling on up to 20,480 cores of the HECToR supercomputer. Overall our work demonstrates that bespoke multigrid preconditioners are highly efficient for solving the pressure correction equation encountered in NWP models.

There are several ways to further improve this work:
so far all tests have been carried out without any orography. It is known that steep gradients can lead to deteriorating performance of the non-linear iteration and we plan to study this by looking at the full non-linear solve for more realistic model problems.
For simplicity we used a finite volume discretisation, but more advanced approaches such as higher-order mixed finite elements can also be used in this framework. This will require the solution of a suitable pressure correction equation in higher order FEM spaces. The parallel performance can also be further improved by, for example, overlapping calculations and communications and strong scaling tests should also be carried out.
Finally, the performance gains from approximate factorisations of the matrix are expected to be significantly higher on GPU systems and hence on such architectures its use may be more justified and more efficient for a wider class of profiles.
\ifpreprint
\section*{Acknowledgements}
\else
\\
\acks
\fi 
This work was funded as part of the NERC project on Next Generation Weather and Climate Prediction (NGWCP), grant numbers NE/J005576/1 and NE/K006754/1. We gratefully acknowledge input from discussions with our collaborators in the Met Office Dynamics Research group and the GungHo! project, in particular Tom Allen, Terry Davies, Markus Gross and Nigel Wood. Ian Boutle kindly provided the aquaplanet UM output used in Section \ref{sec:ResultsAquaplanet}. We would like to thank all DUNE developers and in particular Oliver Sander for his help with extending the parallel scalability of the UG grid implementation.
This work made use of the facilities of HECToR, the UK's national high-per\-for\-mance computing service, which is provided by UoE HPCx Ltd at the University of Edinburgh, Cray Inc and NAG Ltd, and funded by the Office of Science and Technology through EPSRC's High End Computing Programme.

\ifpreprint 
\appendix
\section{Discretisation}\label{sec:Discretisation}
To use a finite volume discretisation we assume that the horizontal grid is divided into elements $T$ of area $|T|$ where the centre of each cell is denoted by $\rhat_T$. The vertical grid is defined by grid levels
\begin{equation*}
  1 = r_0 < r_1 < \dots < r_{n_z} = 1 + H.
\end{equation*}
Each three dimensional grid cell $E$ is defined by a cell $T$ of the horizontal grid and a vertical index $k$, such that cell $E \equiv (T,k)$ with $k=0,\dots,n_z-1$ is bound by $r_k$ and $r_{k+1}$
and we write $r_{k+\frac{1}{2}} = \frac{1}{2}\left(r_{k+1}+r_k\right)$. Any continuous field $u(r,\rhat)$ can be approximated by its average value $u_{T,k}$ in a grid cell as
\begin{equation*}
  \cellint u(r,\rhat) = \cellinthoriz\cellintvert u(r,\rhat)
  \equiv |T| v_k u_{T,k}
\end{equation*}
where we have used
\begin{equation*}
  |E| \equiv \cellint = \cellinthoriz \cellintvert =
  |T| \frac{1}{3}\left(r_{k+1}^3-r_k^3\right) \equiv |T| v_k.
\end{equation*}
As all degrees of freedom in one vertical column are stored consecutively in memory, we also introduce the vector $\mathbf{u}_T$ with $(\mathbf{u}_T)_k = u_{T,k}$.
Then it is straightforward to write down the finite volume discretisation of the individual terms in (\ref{eqn:HelmholtzVector}).
\paragraph{Mass term}
\begin{equation*}
  \cellint \beta u = \cellinthoriz \cellintvert \beta u
  \approx |T| v_k \beta_{T,k} u_{T,k}
\end{equation*}
\paragraph{Diffusion}
\begin{equation}
  \cellint \vec{\nabla}\cdot\left(\mat{\alpha}\vec{\nabla}u\right)
    = \cellinthoriz\cellintvert \left(\frac{1}{r^2}
    \nablatwod\cdot\left(\alpha_{\hS}\nablatwod u\right)
    + \frac{1}{r^2}\dr\left(r^2 \alpha_{r} \dr u\right)
    \right)\label{eqn:discrDiffusion}
\end{equation}
After integration by part in the horizontal direction the first term becomes
\begin{equation*}
  \begin{aligned}
  &\int_{r_{k}}^{r_{k+1}}dr\; \int_{\partial{T}} d\ell\; \alpha_\hS(r,\rhat(\ell)) \vec{n}_{\partial T}\cdot \nablatwod u(r,\rhat(\ell)) \\ &\approx
  (r_{k+1}-r_k) \sum_{T'\in\nbT} \alpha_\hS(r_{k+\frac{1}{2}},\vec{\omega_{TT'}})\frac{|S_{TT'}|\vec{n}_{TT'}\cdot(\rhat_{T'}-\rhat_T)}{|\rhat_{T'}-\rhat_T|^2} \left(u_{T',k}-u_{T,k}\right) \\
  &\equiv \sum_{T'\in\nbT} \left\{(r_{k+1}-r_k)
\frac{|S_{TT'}|\vec{n}_{TT'}\cdot(\rhat_{T'}-\rhat_T)}{|\rhat_{T'}-\rhat_T|^2}
(\alpha_{\hS})_{TT',k}
\right\}\left(u_{T',k}-u_{T,k}\right)
  \end{aligned}
\end{equation*}
where the sum runs over the neighbours $T'$ of the horizontal grid cell. $|S_{TT'}|$ is the length of the edge between the cells $T$ and $T'$ and $\vec{n}_{TT'}$ is an outward normal vector on this edge and tangential to the sphere.
Note that $(\alpha_{\hS})_{TT',k}$ is a field that ``lives'' on the vertical faces of a three dimensional grid cell. The second term in (\ref{eqn:discrDiffusion}) is treated similarly by a vertical integration by parts to obtain
\begin{equation}
  \begin{aligned}
  & |T| \left(
    r_{k+1}^2 \alpha_{r}(r_{k+1},\rhat_T) \sigma_{k+1}
    \frac{u_{T,k+1}-u_{T,k}}{\frac{1}{2}(r_{k+2}-r_{k})}
    -
    r_{k}^2 \alpha_r(r_{k},\rhat_T) \sigma_{k}
    \frac{u_{T,k}-u_{T,k-1}}{\frac{1}{2}(r_{k+1}-r_{k-1})}
  \right)\\
  & \equiv \left\{\sigma_{k+1}|T|\frac{2r_{k+1}^2}{r_{k+2}-r_{k}}(\alpha_{r})_{T,k+1}\right\} \left(u_{T,k+1}-u_{T,k}\right)
   +\left\{\sigma_{k}|T|\frac{2r_{k}^2}{r_{k+1}-r_{k-1}}(\alpha_{r})_{T.k}\right\}\left(u_{T,k}-u_{T,k-1}\right).
  \end{aligned}
  \label{eqn:VerticalDiffusion}
\end{equation}
Often the vertical boundary conditions are fixed by requiring that the vertical velocity is zero at the top and bottom of the atmosphere. This implies mixed boundary conditions of the form $C_1u+C_2\partial_r u=0$ for the pressure.
To simplify the discussion, we use homogeneous Neumann boundary conditions, i.e. $\partial_r u=0$. Because all degrees of freedom in a vertical column are relaxed simultaneously in our solver this should not have any impact on the performance. Furthermore, due to the presence of the zero-order term in (\ref{eqn:HelmholtzVector}) both the individual tridiagonal systems in one column and the elliptic operator are non-singular. In (\ref{eqn:VerticalDiffusion}) Neumann boundary conditions are enforced by setting $\sigma_0 = \sigma_{n_z} = 0$ and $\sigma_k=1$ otherwise. The field $(\alpha_r)_{T,k}$ is defined on the horizontal faces of each three dimensional grid cell.
\paragraph{Vertical advection}
\begin{equation*}
  \begin{aligned}
  \cellint \vec{\xi}\cdot (\vec{\nabla}u) &=
  \cellinthoriz\cellintvert \xi_{r}(r,\rhat) (\dr u(r,\rhat)) \\
  &\approx \frac{1}{2}|T| (r_{k+1}-r_{k})\Big(
    r_{k+1}^2 \xi_{r}(r_{k+\frac{1}{2}},\rhat_T) \frac{u_{T,k+1}-u_{T,k}}{\frac{1}{2}(r_{k+2}-r_{k})}
    + r_{k}^2 \xi_{r}(r_{k-\frac{1}{2}},\rhat_T) \frac{u_{T,k}-u_{T,k-1}}{\frac{1}{2}(r_{k+1}-r_{k-1})}
  \Big) \\
  &= \left\{\sigma_{k+1} |T|\frac{r_{k+1}^2(r_{k+1}-r_k)}{r_{k+2}-r_k}(\xi_{r})_{T,k+1}\right\}\left(u_{T,k+1} - u_{T,k}\right)
  + \left\{\sigma_k|T|\frac{r_{k}^2(r_{k+1}-r_k)}{r_{k+1}-r_{k-1}}(\xi_{r})_{T,k}\right\}\left(u_{T,k} - u_{T,k-1}\right)
  \end{aligned}
\end{equation*}
\subsection{Tridiagonal system}
\label{sec:TridiagonalSystem}
Based on the results in the previous sections we we can now give explicit expressions for the quantities which are needed in (\ref{eqn:TridiagonalNonFactorising}) and (\ref{eqn:TridiagonalFactorising}) to construct the entries of the (tri-) diagonal matrices in (\ref{eqn:TridiagonalSystem}). In the most general case all profile functions depend on the radial coordinate $r$ and the horizontal coordinate $\rhat$. In this case define
\begin{equation}
  \begin{aligned}
    \hat{\beta}_{T,k} &\equiv |T|v_k \beta_{T,k} \\[1ex]
    (\hat{\alpha}_{\hS})_{TT',k} &\equiv \omega^2(r_{k+1}-r_{k}) \frac{|S_{TT'}|\vec{n}_{TT'}\cdot(\rhat_{T'}-\rhat_T)}{|\rhat_{T'}-\rhat_T|^2}(\alpha_{\hS})_{TT',k},&
    (\hat{\alpha}_\hS)_{T,k} &\equiv
    \omega^2\sum_{T'\in\nbT}(\hat{\alpha}_\hS)_{TT',k}\\[1ex]
    (\hat{\alpha}_r)_{T,k} &\equiv \omega^2 \sigma_k|T| \frac{2r_k^2}{r_{k+1}-r_{k-1}} (\alpha_r)_{T,k}, &
    (\hat{\xi}_{r})_{T,k} &\equiv \omega^2\sigma_k |T|\frac{r_k^2(r_{k+1}-r_{k})}{r_{k+1}-r_{k-1}}(\xi_{r})_{T,k}
  \end{aligned}
  \label{eqn:hatFunctionsNonFac}
\end{equation}
which can all be precomputed. If the profile functions can be written as
\begin{xalignat*}{4}
  \alpha_\hS(r,\rhat) &= \alpha_\hS^r(r)\alpha_\hS^\hS(\rhat), &
  \alpha_r(r,\rhat) &= \alpha_r^r(r)\alpha_r^\hS(\rhat), &
  \xi_r(r,\rhat) &= \xi_r^r(r)\xi_r^\hS(\rhat), &
  \beta(r,\rhat) &= \beta^r(r)\beta^\hS(\rhat)
\end{xalignat*}
this implies
\begin{xalignat*}{4}
  (\alpha_\hS)_{T,k} &= (\alpha_\hS^r)_k(\alpha_\hS^\hS)_{T}, &
  (\alpha_r)_{T,k} &= (\alpha_r^r)_k(\alpha_r^\hS)_{T}, &
  (\xi_r)_{T,k} &= (\xi_r^r)_k(\xi_r^\hS)_{T}, &
  \beta_{T,k} &= \beta^r_k\beta^\hS_{T}.
\end{xalignat*}
In this case define
\begin{equation}
  \begin{aligned}
    \hat{\beta}^r_k  &\equiv v_k \beta^r_k, &
    \hat{\beta}^\hS_T   &\equiv |T|\beta^\hS_T \\[1ex]
    (\hat{\alpha}_\hS^r)_k &\equiv (r_{k+1}-r_{k}) (\alpha_\hS^r)_k, &
    (\hat{\alpha}_\hS^\hS)_{TT'} &\equiv \omega^2 \frac{|S_{TT'}|\vec{n}_{TT'}\cdot(\rhat_{T'}-\rhat_T)}{|\rhat_{T'}-\rhat_T|^2}(\alpha_{\hS}^{\hS})_{TT'}\\[1ex]
    & & (\hat{\alpha}_\hS^\hS)_{T} &\equiv
    \sum_{T'\in\nbT}(\hat{\alpha}_\hS^\hS)_{TT'}\\[1ex]
    (\hat{\alpha}_r^r)_k &\equiv \sigma_k \frac{2r_k^2}{r_{k+1}-r_{k-1}} (\alpha_r^r)_k, &
    (\hat{\alpha}_r^\hS)_T &\equiv \omega^2|T| (\alpha_r^\hS)_T \\[1ex]
    (\hat{\xi}_r^r)_k &\equiv \sigma_k \frac{r_k^2(r_{k+1}-r_{k})}{r_{k+1}-r_{k-1}}(\xi_r^r)_k, &
    (\hat{\xi}_r^\hS)_T &\equiv \omega^2|T|(\xi_r^\hS)_T
  \end{aligned}
  \label{eqn:hatFunctionsFac}
\end{equation}
\fi 


\ifpreprint 
\bibliographystyle{unsrt}
\else 
\bibliographystyle{wileyj}
\fi 

\end{document}